\input amstex
\documentstyle{amsppt}
\magnification=\magstep1

\NoBlackBoxes
\TagsAsMath

\pagewidth{6.5truein}
\pageheight{9.0truein}

\long\def\ignore#1\endignore{#1}

\ignore
\input xy \xyoption{matrix} \xyoption{arrow}
          \xyoption{curve}  \xyoption{frame}
\def\edge{\ar@{-}}
\def\dttdar{\ar@{.>}}

\def\dashedge{\ar@{--}}
\endignore

\def\seq{\mathrel{\widehat{=}}}
\def\la{{\Lambda}}
\def\lamod{\Lambda\text{-}\roman{mod}}
\def \len{\operatorname{length}} 
\def\AA{{\Bbb A}}
\def\CC{{\Bbb C}}
\def\PP{{\Bbb P}}
\def\SS{{\Bbb S}}
\def\ZZ{{\Bbb Z}}
\def\NN{{\Bbb N}}
\def\RR{{\Bbb R}}
\def\ZZ{{\Bbb Z}}
\def\hom{\operatorname{Hom}}
\def\soc{\operatorname{soc}}

\def\stab{\operatorname{Stab}}

\def\GL{\operatorname{GL}}

\def\K{\operatorname{K}}

\def\Add{\operatorname{Add}}
\def\Irr{\operatorname{Irr}}
\def\Ext{\operatorname{Ext}}

\def\A{{\Cal A}}
\def\B{{\Cal B}} 
\def\C{{\Cal C}}
\def\D{{\Cal D}}

\def\H{{\Cal H}}
\def\I{{\Cal I}}

\def\T{{\Cal T}}

\def\M{{\Cal M}}
\def\fancyO{{\Cal O}}
\def\P{{\Cal P}}

\def\R{{\Cal R}}
\def\S{{\sigma}}
\def\s{{\frak s}}
\def\T{{\Cal T}}
\def\U{{\Cal U}}

\def\e{\bold {e}}
\def\Phat{\widehat{P}}
\def\aut{\operatorname{Aut}}
\def\Aut{\operatorname{Aut}}
\def\autlap{\operatorname{Aut}_\la(P)}
\def\End{\operatorname{End}}

\def\K{\operatorname{K}}

\def\dim{\operatorname{dim}}
\def\locdim{\operatorname{loc\,dim}}

\def\Schu{\operatorname{Schu}}

\def\Aff{\operatorname{Aff}}
\def\AffS{\Aff(\S)}

\def\modlad{\operatorname{\bold{Mod}}^\Lambda_d}
\def\toptd{\operatorname{\bold{Mod}}^T_d}

\def\Hom{\operatorname{Hom}}
\def\grasstd{\operatorname{\frak{Grass}}^T_d}
\def\grassS{\operatorname{\frak{Grass}}(\S)}
\def\grassSS{\operatorname{\frak{Grass}}(\SS)}
\def\grass#1{\operatorname{\frak{Grass}}(#1)}
\def\Gr{\operatorname{Gr}}
\def\Var{\operatorname{Var}}

\def\degen{\le_{\text{deg}}}
\def\underbardim{\operatorname{\underline{dim}}}
\def\id{\operatorname{id}}
\def\term{\operatorname{end}}
\def\start{\operatorname{start}}
\def\vsubseteq{\cup\kern1.0pt{\vrule width0.35pt height1.38ex}}

\def\irred{{\bf 1}}
\def\BongAdv{{\bf 2}}
\def\Bong{{\bf 3}}
\def\Bongtrond{{\bf 4}}
\def\GeomII{{\bf 5}}
\def\GeomIV{{\bf 6}}
\def\Bor{{\bf 7}}
\def\Hil{{\bf 8}}
\def\GeomI{{\bf 9}}
\def\preview{{\bf 10}}
\def\Hum{{\bf 11}}
\def\Ki{{\bf 12}}
\def\Kra{{\bf 13}}
\def\BruScho{{\bf 14}}
\def\GIT{{\bf 15}}
\def\New{{\bf 16}}
\def\Rie{{\bf 17}}
\def\RosII{{\bf 18}}
\def\Ros{{\bf 19}}
\def\Scho{{\bf 20}}
\def\Zwa{{\bf 21}}

\topmatter

\title  Classifying representations by way of
Grassmannians
\endtitle

\rightheadtext{Classifying finite dimensional representations}

\author Birge Huisgen-Zimmermann \endauthor

\address Department of Mathematics, University of California, Santa
Barbara, CA 93106 \endaddress

\email birge\@math.ucsb.edu \endemail

\thanks This research was partially supported by a grant from the
National Science Foundation. \endthanks

\dedicatory Dedicated to the memory of Sheila Brenner \enddedicatory

\abstract Let $\la$ be a finite dimensional algebra over an
algebraically closed field.  Criteria are given which characterize
existence of a fine or coarse moduli space classifying, up to
isomorphism, the representations of
$\la$ with fixed dimension $d$ and fixed squarefree top $T$.  Next to
providing a complete theoretical picture, some of these equivalent
conditions are readily checkable from quiver and relations.  In case
of existence of a moduli space  --  unexpectedly frequent in light of
the stringency of fine classification  --  this space is always
projective and, in fact, arises as a closed subvariety
$\grasstd$ of a classical Grassmannian.  Even when the full moduli
problem fails to be solvable, the variety $\grasstd$ is seen to have
distinctive properties recommending it as a substitute for a moduli
space.  As an application, a characterization of the algebras having
only finitely many representations with fixed simple top is obtained; 
in this case of `finite local representation type at a given simple
$T$', the radical layering $\bigl( J^lM/ J^{l+1}M \bigr)_{l \ge 0}$ is
shown to be a classifying invariant for the modules with top $T$.  This
relies on the following general fact obtained as a byproduct: Proper
degenerations of a local module $M$ never have the same radical
layering as $M$. 
\endabstract

\endtopmatter

\document

\head 1. Introduction and terminology \endhead

Throughout, $\la$ denotes a
finite dimensional algebra over an algebraically closed field $K$. 
While the full representation-theoretic picture of $\la$ is beyond the
scope of a complete description if $\la$ has wild representation type,
substantial portions of this representation theory do lie within reach
in many wild situations.  The part we address here consists of the
representations $M$ of fixed dimension
$d$ that have fixed squarefree top
$T = M/JM$, where $J$ is the Jacobson radical of $\la$; in other
words, we require that the simple (left)
$\la$-modules occur with multiplicity at most $1$ in $T$.  Our
primary goal is to decide when the restricted classification
problem for isomorphism classes of
$d$-dimensional representations with top $T$ has a coarse or fine moduli
space.  In rough terms, this means that we seek to bijectively
parametrize these isomorphism classes by the points of an algebraic
variety such that the structure constants `evolve Zariski continuously'
as one moves along the parameter space   --  this continuity condition
is made precise in the concept of a family  --  and so that other
continuous parametrizations are uniquely induced by the distinguished
one  --  the latter requirement is made precise via universal
properties of varying degrees of stringency.  (This classification
philosophy actually predates the rigorous definition of coarse and fine
moduli spaces given by Mumford in the 1960s; it already underlies
Riemann's classification of nonsingular projective curves in the 1850s,
where the term `moduli' was coined, standing for a
`structure-determining collection of continuous parameters in $\CC$'.) 
A precise definition of a fine moduli space can be found at the end of
Section 1.

In our situation, the existence problems for fine and coarse moduli
spaces will turn out to be equivalent.  One of our main objectives is
to provide readily verifiable necessary and sufficient conditions for
existence.  In the positive case, we exhibit the fine
moduli space --  it is always projective  --  together with the
universal family that classifies, up to isomorphism, the
$d$-dimensional modules with top $T$.
\medskip

In the classical affine variety $\modlad$ whose $\GL_d$-orbits
bijectively parametrize the isomorphism classes of $d$-dimensional
$\la$-modules, all orbits corresponding to non-semisimple modules fail
to be closed.  Therefore, the standard methods of invariant theory,
typically restricting attention to the closed orbits, are a priori not
helpful.  In a seminal article,
\cite{\Ki}, King coped with this difficulty by adapting  Mumford's
concept of stability of vector bundles on projective curves, which led
him to focus on
$\la$-modules which are `(semi-)stable' relative to a given additive
function $\Theta: K_0(\lamod) \rightarrow \RR$.  For the
$\Theta$-stable modules, a fine moduli space classifying up to
isomorphism is guaranteed, but this class of representations is
often hard to identify or assess in size.  

Here we initiate a
quite different approach to the moduli problem.  Let us start by
presenting a projective variety whose points always parametrize, not
necessarily bijectively, the
$d$-dimensional representations with fixed top $T$.  As in the
classical setting, the isomorphism classes are again in 1-1
correspondence with the orbits of an algebraic group action, but these
orbits have lower dimension and  are closed much more frequently.  This
variety was first introduced by Bongartz and the author in
\cite{\GeomII} and
\cite{\GeomIV}:  Fix
$d$ and a semisimple module $T$, not yet assumed to be squarefree, with
projective cover
$P$. We denote by
$\grasstd$ the closed subvariety of the classical Grassmannian
of all
$(\dim P - d)$-dimensional subspaces of $JP$ that consists of the
$\la$-submodules of $JP$ of dimension $\dim P - d$.  Clearly,
$\grasstd$ is a projective variety endowed with a morphic
action of the automorphism group $\Aut_{\la}(P)$, the orbits of which
coincide with the fibres of the surjection
$$\phi: \grasstd \longrightarrow \{ \text{isomorphism types of\ }
d\text{-dimensional modules with top\ } T\}$$ 
sending $C$ to the isomorphism type
of $P/C$.  See Section 2 for basic facts about the
variety $\grasstd$ and its relationship to its classical counterpart.
\medskip 

From now on, we assume $T$ to be squarefree.  Without this
hypothesis, the solution to the moduli problem requires modification and
an additional layer of machinery. It will be treated in a sequel to the
present work; a brief summary is given in Section 7. 
Adopting the notion of a `family of
$\la$-modules' introduced by King (cf\.
\cite{\Ki} and the last paragraph of Section 1), we obtain the
following characterization of the triples $\la$,
$T$,
$d$ for which the moduli problem has a solution (see Theorems 4.2 and
4.4 for somewhat stronger results and Section 2\.A for the concept of a
degeneration):

\proclaim{Theorem A}  The following statements are equivalent for
squarefree $T$.
\roster
\item There exists a {\it coarse} moduli space classifying the
$d$-dimensional $\la$-modules with top $T$, up to isomorphism.

\item  There exists a {\it fine} moduli space classifying the
$d$-dimensional $\la$-modules with top $T$, up to isomorphism.

\item No $d$-dimensional $\la$-module with top
$T$ has a proper degeneration with top $T$.

\item Every
submodule $C \subseteq JP$ of codimension $d$ in
$P$ is invariant under all endomorphisms of $P$. 

\item $\grasstd$ {\it is} the fine moduli space for the
$d$-dimensional $\la$-modules with top $T$ {\rm{(}}in particular, all
$\autlap$-orbits of $\grasstd$ are singletons{\rm{)}}.
\endroster
\endproclaim

We emphasize that Theorem A addresses the existence of a moduli space
for the isomorphism classes of {\it all\/} $d$-dimensional
$\la$-modules with top $T$; in particular, no stability conditions are
imposed.

The equivalence of (3) and (4) holds `pointwise', in the
following sense: An individual module
$M$ with top $T$, say $M = P/C$, is devoid of proper top-$T$
degenerations precisely when $C$ is fully invariant in
$P$; the former condition is, a priori, hard to check, whereas
recognizing full invariance is easy.  Invariance is a strong requirement,
but nonetheless globally satisfied for the submodules of $JP$
of fixed codimension in a wide range of interesting cases  -- 
see Corollary 4.5 for illustration. (For the sake of
contrast:   If
$d \ge 2$, the full collection of $d$-dimensional modules has a coarse
moduli space only when $\la$ is semisimple.)   

When $T$ is simple and the Gabriel quiver of $\la$ has no oriented
cycles, the existence of fine moduli spaces for representations with
fixed top $T$ can alternatively be deduced from King's method; for
detail, see the comments following Theorem 4.4.   Moreover, we refer to
work  of LeBruyn and Schofield, \cite{\BruScho} and
\cite{\Scho}, for results concerning the structure of
moduli spaces over hereditary algebras in case of existence;  in
particular, Schofield provides very general sufficient conditions for
rationality in that case \cite{\Scho}. 

One of the representation-theoretically most useful features of
$\grasstd$  --  here and in subsequent work 
--  is a finite affine cover, $\bigl(\grassS
\bigr)_\sigma$, with the following properties (see Theorem 3.5): (a) It
consists of $\autlap$-stable charts, and (b), it is the disjoint union
of {\it open\/} affine covers of those subvarieties $\grassSS$ of
$\grasstd$ which consist of the points representing modules $M$ with
fixed radical layering
$\SS = \bigl(J^l M/J^{l+1}M \bigr)_{l \ge 0}$; the latter is meaningful
as we identify isomorphic semisimple modules.  For a bit more detail,
suppose $\la = KQ/I$, where
$Q$ is a quiver and $I$ an admissible ideal in the path algebra $KQ$. 
In this scenario, each chart $\grassS$ is indexed by a set $\S$ of
paths tied to the algebraic structure of the $\la$-modules
corresponding to the points of
$\grassS$; in fact, given a module $M$, the family of all eligible sets
$\S$ is an isomorphism invariant of $M$.   The cover $\bigl( \grassS
\bigr)_\sigma$ is  distinguished by its functoriality relative to the
ideal
$I$ of relations of $\la$, the pertinent functor depending only on $Q$,
$d$, and the set of vertices determining
$T$ (3.17).  We add a few
comments relating the $\grassS$ to the Schubert
cells of the  encompassing full Grassmannian
$\operatorname{Gr}(\dim P - d, JP)$: The intersections of
$\grasstd$ with Schubert cells usually fail to satisfy the above
conditions (a) and (b).  In
particular, the intersections with open Schubert
cells are hardly ever stable under the $\autlap$-action.  On the
other hand, each chart $\grassS$
embeds into a suitable open Schubert cell of $\operatorname{Gr}(\dim P
- d, JP)$, and the affine coordinates of
$\grassS$ introduced in 3.10 are essentially the Pl\"ucker coordinates
that correspond to this embedding.  However, the Pl\"ucker coordinates
can be pared down to a comparatively small subfamily
encoding all relevant information, due to our specific algebraic
setting; it is this economy that permits comparatively effortless
analysis of even large examples (see \cite{\preview}).   While the cover
$\bigl(\grassS\bigr)_\S$ provides the backbone for our proofs,
familiarity with its specifics is not required for an understanding of
the main results.  The reader only interested in the theorems and their
theoretical applications is therefore encouraged to skip the somewhat
technical portion of 3.B following Theorem 3.5.

Polynomials defining any of the affine subvarieties $\grassS$ of
$\grasstd$ are available from quiver and relations of $\la$ by way of
easy combinatorial manipulations (3.13).  These polynomials provide
the foundation for all of our concrete examples here and in
\cite{\preview, \irred}; in particular, they make all of
the geometric conditions arising in the theorems checkable through
Gr\"obner methods (see \cite{\irred}).  Furthermore, they
establish a combinatorial link between the geometry of the affine
charts $\grassS$ and algebraic features of the classes of
representations parametrized by them (see, e.g.,
\cite{\preview}, Section 5).  
          
In situations where the $d$-dimensional $\la$-modules with fixed top
fail to have a moduli space, the natural next step is to  subdivide the
target class by fixing further discrete invariants.  There is an obvious
refinement of our primary partition of the $d$-dimensional
representations $M$ in terms of tops:  namely, the partition in terms of
full radical layerings $\SS = \SS(M) = (J^l M/J^{l+1}M)_{0
\le l \le L}$, where  $L$ is chosen so that $J^{L+1} = 0$; in other
words, we fix the  matrix recording the multiplicities of {\it all}
simple composition factors of
$M$ in a format that keeps  track of their `layer-locations'. The main
benefit of restricting one's focus to
$\grassSS$ lies in the fact that the relative closures of the
$\autlap$-orbits in $\grassSS$ are better understood than the
closures in $\grasstd$ (see Section 4); in particular, the orbits of 
$\grassSS$ under the action of the unipotent radical of $\autlap$ are
always relatively closed.  When specialized to the case of a simple
top  $T$, this yields:  

\proclaim{Theorem B} {\rm(}See Section {\rm 2.A} for
terminology\/{\rm).}  No $\la$-module $M$ with simple top has a  proper
degeneration $M'$ with
$J^l M/ J^{l+1}M \cong J^l M'/ J^{l+1}M'$ for $0 \le l \le L$.
\endproclaim 

Section 4 also displays a first set of reasons for
viewing $\grasstd$, and a fortiori $\grassSS$, as a useful approximation
to a moduli space for the corresponding classification problem whenever
the moduli problem in Mumford's sense is not solvable (for further
backing of this viewpoint, we refer to
\cite{\preview}).  To name one reason, the fibres $\autlap.C$ of the
map $\phi$ are well-understood in terms of their intrinsic structure:  
If $T$ is simple, the $\autlap$-orbits are affine spaces $\AA^m$,
where $m$ is bounded above by the multiplicity of $T$ in
$JP$.  For general squarefree $T$, each orbit
$\autlap.C$ is a direct product of an affine space $\AA^m$ and a
torus $(K^*)^r$, where $\dim T - r$ is the number of indecomposable
summands of $P/C$. 

Our characterization of the algebras that have finite
local representation type with respect to a simple module
$T$, as defined in the abstract, hinges on Theorem B.  It
lends additional support to the central role we attribute to the radical
layering as a discrete invariant of a representation.  Condition (3)
below is  decidable from quiver and relations of $\la$ by way of the
mentioned polynomials; if it is satisfied, the modules with top $T$ can
be explicitly constructed from these data.  See Theorem 5.2 for
additional information.  

\proclaim{Theorem C}  For
any simple
$\la$-module
$T$, the following conditions are equivalent {\rm{(}}they are not
left-right symmetric  -- we refer to left modules, say\/{\rm{)}}:

\roster
\item There are only finitely many $\la$-modules with top $T$, up
to isomorphism.

\item  If $M$ and $N$ are $\la$-modules with top $T$, then $M \cong N$
if and only if $\SS(M) = \SS(N)$.

\item For every $d\in\NN$ and every sequence $\SS = (\SS_0, \SS_1,
\dots, \SS_L
\bigr)$ of semisimple $\la$-modules with $\SS_0 = T$ and $\sum_i \dim
\SS_i = d$, the subset $\grassSS$ of $\grasstd$ is either empty or
irreducible of dimension 
$$\bigl(\text{\rm  multiplicity of}\ T \ \text{\rm in\ } M \bigr)\  -\
\dim \End_\la(M)$$  
for some {\rm{(}}equivalently, for all\/{\rm{)}}
modules
$M$ with radical layering $\SS$.
\endroster
\endproclaim

While finite local representation type at $T$ forces $\SS(-)$ to
separate isomorphism classes, the dimension vector fails to separate in
general (Example 5.4).

\medskip

{\bf Notation and terminology.}  Throughout, we will assume $\la$ to
be basic.  Due to algebraic closedness of the ground field $K$, we
may thus, without loss of generality, assume that
$\la = KQ/I$, where $Q$ is a quiver, and
$I$ an admissible ideal in the path algebra $KQ$; the latter means that
$I$ is contained in the ideal generated by all paths of length $2$ and
contains some power of this ideal.  The quiver provides us with a
convenient set of primitive idempotents
$e_1, \dots, e_n$ of
$\la$, which are in bijective correspondence with the vertices of
$Q$; we will, in fact, not distinguish between the vertices and the
$e_i$.  As is well-known, the factors $S_i =
\la e_i /J e_i$ form a set of representatives of the simple (left)
$\la$-modules.  By
$L$ we denote the largest integer for which the power
$J^L$ of the Jacobson radical does not vanish; in other words, $L+1$ is
the Loewy length of
$\la$. 

 Given any (left) $\la$-module $M$, an element $x \in M$ will be called a
{\it top element\/} of
$M$ if $x \notin JM$ and $x$ is normed by some $e_i$, that is $x = e_i x$
for some $i \in \{1, \dots, n\}$.  The isomorphism invariant 
$$\SS(M) = (M/JM,\, JM/J^2M,\, \dots,\, J^{L-1}M/J^L M,\, J^L M)$$ of
$M$ will be referred to as the {\it sequence of radical layers of\/}
$M$, or, more briefly, the {\it radical layering\/} of $M$.  

Moreover, we observe the following conventions:  The product $pq$ of two
paths $p$ and $q$ in $KQ$ stands for `first
$q$, then $p$' if $\term (q) = \start (p)$, and zero otherwise  (so, in
particular, $p = pe_i$ means that the path
$p$ starts in the vertex $e_i$).  In line with this notation, we call a
path
$p_1$ a {\it right subpath\/} of $p$ if $p = p_2 p_1$ for some path
$p_2$.  We will generally gloss over the distinction between the left
$\la$-structure of $M \in \lamod$ and the induced left
$KQ$-structure; in particular, we let paths operate on $\la$-modules
without using residue notation.

For background on moduli problems, we refer to
\cite{\New}, but recall the definition of a fine moduli space for our
specific problem.  Our concept of a family of $\la$-modules
is that introduced by King in
\cite{\Ki}:  Namely, a {\it family of $d$-dimensional $\la$-modules 
parametrized by an algebraic variety
$X$\/} is a pair $(\Delta, \delta)$, where $\Delta$ is a (geometric)
vector bundle of rank $d$ over $X$ and $\delta: \la \rightarrow
\End(\Delta)$ a
$K$-algebra homomorphism.  Our notion of {\it equivalence of families\/}
parametrized by the same variety 
$X$, finer than King's in general, is the coarsest possible to separate
isomorphism classes:  Namely,
$(\Delta_1, \delta_1) \sim (\Delta_2, \delta_2)$ precisely when, for
each $x \in X$, the fibre of $\Delta_1$ over $x$ is $\la$-isomorphic to
the fibre of $\Delta_2$ over $x$.  As is common, given a family
$(\Delta, \delta)$ parametrized by
$X$ and a morphism $\tau: Y
\rightarrow X$ of varieties, the {\it induced family\/} $\tau^*(\Delta,
\delta)$ over $Y$ is the pullback of $(\Delta, \delta)$ along $\tau$. 
In this context, a variety
$X$ is a {\it fine moduli space\/} for (families of) $d$-dimensional
modules with top $T$ if there exists a family $(\Gamma, \gamma)$ of such
modules parametrized by $X$ which has the property that an arbitrary
family  -- parametrized by $Y$ say  --  is equivalent to a family
$\tau^*(\Gamma, \gamma)$ induced via a unique morphism $\tau: Y
\rightarrow X$;  accordingly $\Gamma$ is called the {\it universal\/}
family in case of existence.  In particular, the requirements on
$\Gamma$ entail that every
$d$-dimensional module with top $T$ is isomorphic to precisely one
fibre of the bundle
$\Gamma$.  For the more common definition of a fine moduli space through
representability of a suitable functor, as well as for the concept of  a
{\it coarse moduli space\/}, see \cite{\New, pp\. 23, 24}.    
\bigskip

\head 2. The varieties $\toptd$ and $\grasstd$ for general $T$ \endhead

We fix $d \ge 1$ and a finitely generated semisimple module $T$, say $T =
\bigoplus_{ 1 \le i \le n} (\la e_i/ Je_i)^{t_i}$.  Moreover, we let
$P=\bigoplus_{1 \le i \le n} (\la e_i)^{t_i}$ be the projective cover of
$T$.  In the first three subsections, we discuss the classical and
Grassmannian varieties associated to the $d$-dimensional modules with
top $T$, along with connections between them.

\subhead 2.A. The classical setup \endsubhead

If we denote by $Q^*$ the union of the set $\{e_1, \dots, e_n\}$ of
vertices with the set of arrows of
$Q$, then $Q^*$ generates $\la$ as a $K$-algebra, and the traditional
variety of
$d$-dimensional left
$\la$-modules, $\modlad$, takes on the form   
$$\modlad = \{(x_\alpha)_{\alpha \in Q^*} \in \prod _{\alpha \in Q^*}
\operatorname{End}_K(K^d) \mid \ \text{the\ } x_\alpha \
\text{satisfy all the rel'ns of the $\alpha$ in\ } \la\}.$$ This affine
algebraic variety carries a morphic $\GL_d$-action by conjugation
accounting for change-of-basis transformations.  Clearly, the fibres of
the representation map
$$R: \modlad \longrightarrow \{\text{isomorphism types of\ }
d\text{-dimensional\ } \la\text{-modules}\}$$ are precisely the
$\GL_d$-orbits of $\modlad$. The slice of $\modlad$ which will be
pivotal for our investigation is the locally closed subvariety
$$\toptd = {\operatorname{\bold{Mod}}}^{\la,T}_d = \{x \in \modlad \mid
R(x)/JR(x) \cong T  \}.$$  
Clearly, $\toptd$ is stable under the $\GL_d$-action, and its
$\GL_d$-orbits bijectively parametrize the $d$-dimensional modules with
top $T$.  Hence, the plausible next step in looking for a moduli space
is to attempt factoring the $\GL_d$-action out of $\toptd$.

We recall two concepts of
quotient of an algebraic variety modulo an algebraic group action which
will be crucial in the sequel:  Suppose that
$X$ is a variety endowed with a morphic action of an algebraic group
$G$.  Then a {\it categorical quotient} of $X$ by $G$ is a morphism of
varieties
$\pi: X \rightarrow Z$ which is constant on the $G$-orbits of $X$ and
satisfies the following universal property:  Every morphism
$\psi:X \rightarrow Y$ that is constant on the orbits factors uniquely
through $\pi$.  Moreover, a morphism
$\pi: X \rightarrow Z$ having fibres coinciding with the $G$-orbits is
called a {\it geometric quotient\/} of $X$ by $G$ if it is surjective and
open, and satisfies the following condition relating the structure sheaf
of $Z$ to the rings of invariants for the $G$-action: for every open
subset $U$ of
$Z$, the comorphism $\pi^*$ induces an isomorphism from the ring
$\fancyO_Z (U)$ of regular functions on $U$ to the ring
$\fancyO_X (\pi^{-1}(U))^G$ of $G$-invariant regular functions on
$\pi^{-1}(U)$.  Every geometric quotient of $X$ by $G$ is a categorical
quotient and hence unique in case of existence.  It is, moreover,
well-known that, in case $X$ is affine and $G$ reductive (e\.g\., in
case $X = \modlad$ and $G = \GL_d$), a categorical quotient of $X$ by
$G$ always exists; this quotient is geometric precisely when all
$G$-orbits of $X$ are closed.  Sufficiency of this latter condition does
not extend to quasi-affine varieties, however  --  for counter-examples
based on subvarieties of
$\modlad$, see Section 4.

We begin with a well-known criterion characterizing the
existence of a coarse moduli space, specialized to a class of modules
represented by a subvariety (= locally closed subset) of
$\modlad$.  It is due to Mumford \cite{\GIT}, but we refer to
Newstead's formulation in \cite{\New}, since we have adopted the
terminology of that text.

\proclaim{Criterion 2.1}  Let $X$ be a $\GL_d$-stable subvariety of 
$\modlad$ and $\C$ the collection of modules
represented by $X$. Then there exists a coarse moduli space classifying
the modules in $\C$ up to isomorphism if and only if $X$ has a
categorical quotient $\pi: X \rightarrow X/ \GL_d$ which is an orbit
map; the latter means that the map $\pi$ is surjective and that its
fibres coincide with the $\GL_d$-orbits of $X$.  In the positive case,
the quotient
$X / \GL_d$ is the coarse moduli space we are looking for.
\endproclaim

\demo{Proof} By \cite{\New, Proposition 2.13}, it suffices to specify a
family $(\varUpsilon, \varepsilon)$ of
$\la$-modules, parame\-trized by $X$, which has the local universal
property relative to families of
$d$-dimensional modules with top $T$.  Local universality means that, for
an arbitrary family
$(\Delta, \delta)$  -- parametrized by some variety $Y$ say  --  and for
any point $y \in Y$, there exists a neighborhood $U$ of $y$ such that
$\Delta |_U$ is induced from $(\varUpsilon, \varepsilon)$ by way of some
(not necessarily unique) morphism $U \rightarrow
X$.  The obvious candidate for $(\varUpsilon, \varepsilon)$ is as
follows:  Take $\varUpsilon$ to be the trivial bundle of rank $d$ over
$X$, and define $\varepsilon: \la
\rightarrow \End(\varUpsilon)$ by the requirement that, for
$\alpha \in Q^*$, the endomorphism $\varepsilon(\alpha)$ of $\varUpsilon$
coincides with
$x_{\alpha}$ on the fibre above $x$.  To verify the local universal
property, it obviously suffices to show that every family $(\Delta,
\delta)$ based on a trivial bundle $\Delta = U \times K^d$ is induced
from
$(\varUpsilon, \varepsilon)$.  To do so, we let $\tau: U \rightarrow
X$ be the morphism which sends any element $y$ in $U$ to the point
$x \in X$ with the property that, for each
$\alpha \in Q^*$, the endomorphism $x_\alpha$ equals the restriction of
$\delta(\alpha)$ to  the fibre above $y$.  Then $(\Delta, \delta)$ is
equivalent to
$\tau^*(\varUpsilon, \varepsilon)$ as required.  \qed
\enddemo

A  necessary condition for the existence of a coarse moduli space for
the $d$-dimensional representations with top $T$ is thus immediate: 
Namely, all $\GL_d$-orbits of
$\toptd$ need to be closed (in standard jargon, this amounts to
excluding the `jump phenomenon').  In representation-theoretic terms,
this condition says that no
$d$-dimensional module $M$ with top $T$ has a proper degeneration
sharing that top.     

Recall that a {\it degeneration\/} of a
$d$-dimensional module
$R(x)$ is a module $R(y)$ with the property that $y$ belongs to the
closure of the orbit $\GL_d.x$ in $\modlad$.  One writes $R(x)
\degen R(y)$ in that case to reflect the fact that `degenerates to' is a
partial order on the isomorphism classes of $d$-dimensional modules. 
For background on degenerations, we point to \cite{\BongAdv,
\Bongtrond, \Kra, \Rie, \Zwa}.  In view of Criterion 2.1, our present
investigation of moduli spaces automatically involves existence questions
for degenerations.  As byproducts, we will thus obtain some preliminary
results on what we call `top-stable' and `layer-stable' degenerations. 
(This thread will be picked up in
\cite{\preview}.) 

\definition{Definition 2.2}  A degeneration $M'$ of a finitely generated
left $\la$-module $M$ is called {\it top-stable\/} if $M/JM
\cong M'/JM'$, {\it layer-stable\/} in case $J^l M/J^{l+1}M \cong J^l
M'/J^{l+1}M'$ for all $l$.
\enddefinition

Clearly, the layer-stable degenerations of $M$ are a fortiori
top-stable.  Moreover, it is obvious that the top-stable degenerations
of a module $M = R(x)$ with
$x \in \toptd$ are precisely those $R(y)$ for which $y$ belongs to the
relative closure of $\GL_d.x$ in
$\toptd$.  Analogously, if $\SS = \SS(M)$, then the layer-stable
degenerations of $M$ coincide with the $R(y)$ where $y$ traces the
relative closure of $\GL_d.x$ in the subvariety 
$$\operatorname{\bold{Mod}}(\SS) = \{ x \in \toptd \mid R(x)\ \text{has
radical layering}\ 
\SS\}.$$

\subhead 2.B. The Grassmannian setup \endsubhead 

Keep the notation introduced at the beginning of Section 2, and set $d'
=\dim_K\nobreak P\nobreak -\nobreak d$. We define $\grasstd$ to be the
closed subvariety of the Grassmannian of $d'$-dimensional subspaces of
$JP$ consisting of those subspaces that are $\la$-submodules of $JP$.  
In other words,
$\grasstd$ contains precisely those
$\la$-submodules $C$ of  $JP$ for which $P / C$ is a $d$-dimensional
module (automatically having top $T$).  Clearly,
$\grasstd$ is a closed subvariety of the usual Grassmannian of
$d'$-dimensional subspaces of $JP$ and, as such, is in turn projective. 
Note moreover that the variety
$\grasstd$ comes coupled with a surjective map 
$$\phi = \phi^T_d : \grasstd \longrightarrow \{ \text{isom.~types
of\ } d\text{-diml.~modules with top\ } T\}, \quad C \mapsto [P/C].$$  
We will refer to this map as the {\it representation map of the
Grassmannian of $d$-dimensional modules with top $T$\/}.  As already
pointed out, the fibres of $\phi$ coincide with the orbits of the
natural morphic action of $\aut_{\la}(P)$ on
$\grasstd$, denoted by $f.C$.   The fact that the $\aut_{\la}(P)$-orbits
of $\grasstd$ will turn out to be structurally simpler than the
$\GL_d$-orbits of $\toptd$ in many cases, hinges on the following
elementary fact. 

\definition{Observation 2.3}  $\aut_\la (P) \cong \aut_\la^u (P) \rtimes
\aut_\la (T)$, where 
$$\aut_\la^u (P) = \{ \id + f_0 \mid f_0 \in \hom_\la (P,JP)\}$$ is the
unipotent radical of
$\aut_\la (P)$.  The automorphism group $\aut_\la (T)$ of the top $T$ of
$P$ is isomorphic to
$\prod_{1 \le i \le n} {\GL}_{t_i}(K)$, a torus if all the $t_i$ are
at most $1$. 

The isomorphism $\aut_\la (P) \cong \aut_\la^u (P) \rtimes
\aut_\la (T)$ can be replaced by equality if, for instance, we
identify $\aut_\la(T)$ with the subgroup $\H \le \autlap$ consisting of
those $\la$-automorphisms of $P$ which leave the subspace $\bigoplus_{1
\le i \le n} (Ke_i)^{t_i}$ of
$P =  \bigoplus_{1 \le i \le n} (\la e_i)^{t_i}$ invariant.   \qed   
\enddefinition

We use hereditary algebras for two simple illustrations: First let $\la
= KQ$, where $Q$ is the quiver \
\ignore $\xymatrixrowsep{2pc}\xymatrixcolsep{3pc}
\vcenter{\xymatrix{ 2 &1 \ar[l]\ar@/^/[l]\ar@/_/[l]   }}$ \endignore
\ .  Then $\operatorname{\frak{Grass}}^{S_1}_2 \cong
\operatorname{\frak{Grass}}^{S_1}_3 \cong \PP^2$.  

Now suppose that 
$\la = KQ$, where $Q$ is the quiver \
\ignore $\xymatrixrowsep{2pc}\xymatrixcolsep{3pc}
\vcenter{\xymatrix{ 2 &1 \ar[l]\ar@/^/[l]\ar@/_/[l]
\ar[r]\ar@/^/[r]\ar@/_/[r] &3  }}$ \endignore\ .  Then 
$\operatorname{\frak{Grass}}^{S_1}_4$ consists of four irreducible
components (which coincide with the connected components), two of them
copies of $\PP^2 \times \PP^2$ and the other two singletons.   In all of
these instances, the $\autlap$-orbits are reduced to points, that is,
the representation map $\phi$ is a bijection from $\frak{Grass}^{S_1}_4$
to the isomorphism types of $4$-dimensional modules with top
$S_1$.

Yet, in the following example, it is not. 

\definition{Example 2.4}  Let $\la = KQ/\langle
\omega^2 \rangle$, where
$Q$ is the quiver

\ignore
$$\xymatrixrowsep{1pc}\xymatrixcolsep{3pc}
\xymatrix{ 1 \ar@(ul,dl)_{\omega} \ar[r]^{\alpha} &2 }$$
\endignore

\noindent  Then $\frak{Grass}^{S_1}_3 \cong \PP^1$, but
there are only two isomorphism classes of $d$-dimensional modules with
top $S_1$, namely $\la e_1/ \la(\alpha - k \alpha \omega)$ for $k \in K$
and
$\la e_1/ \la \alpha \omega$.
\qed\enddefinition  

In general, it is difficult to determine the structure of
$\grasstd$ `at one blow'.  This is one of the reasons for subdividing 
$\grasstd$ into more accessible affine portions, still stable under the
$\autlap$-action.  Provided that they are nicely located within
$\grasstd$, one can then piece together the representation-theoretic
information encoded in $\grasstd$ from the collection of affine frames. 
Such a subdivision will be introduced in two stages (Sections 2.D and
3.B).   

\subhead 2.C.  Connection between $\toptd$ and $\grasstd$
\endsubhead     

Geometric information on the orbits and orbit closures of the
$\aut_\la(P)$-action of
$\grasstd$ smoothly translates into information about the orbits and
orbit closures of the $\GL_d$-action of the classical variety
$\toptd$, as spelled out in the following proposition.  It is due to
Bongartz and the author \cite{\GeomIV, Proposition C} and was
previously applied towards an investigation of uniserial modules, which
triggered many of the ideas developed here.          

\proclaim{Proposition 2.5}  Consider the map from the set of all
$\autlap$-orbits of $\grasstd$ to the set of $\GL_d$-orbits of
$\toptd$, sending any $\autlap$-orbit $\autlap.C = \phi^{-1}(P/C)$ to the
$\GL_d$-orbit
$R^{-1}(P/C)$.  This map extends to a one-to-one inclusion-respecting
correspondence between the set of
$\aut_{\la}(P)$-stable subsets of
$\grasstd$ on one hand, and the set of
$\GL_d$-stable subsets of $\toptd$ on the other, which preserves
openness, closures, irreducibility, connectedness, and types of
singularities.   

Furthermore:  Suppose that $V$ is an $\autlap$-stable subvariety of
$\grasstd$ and $W$ the corresponding $\GL_d$-stable subvariety of
$\toptd$.  Then: 

{\rm (a)} $V$ has a geometric quotient by its
$\autlap$-action  if and only if $W$  has a geometric quotient by its
$\GL_d$-action, and in the positive case, these quotients coincide.  

{\rm (b)}  $V$ has a categorical quotient modulo $\autlap$ if and only
if $W$ has a categorical quotient modulo $\GL_d$.  In the positive case,
these quotients are isomorphic.  The fibres of the former are equal to
the $\autlap$-orbits of $V$ if and only if the fibres of the latter are
equal to the $\GL_d$-orbits of $W$. \endproclaim

\demo{Proof}  Only the statements under (b) were not addressed in
\cite{\GeomIV}.  They are immediate consequences of Lemmas 1 and 2
preceding \cite{\GeomIV, Proposition C}.  Indeed, a variety $X$ endowed
with an action of an algebraic group $H$ having the following properties
is presented there:  $H$ contains
$\autlap$ and $\GL_d$ as subgroups meeting only in $1$, as well as two
closed normal subgroups $N_1$ and $N_2$ such that $H/N_1 \cong
\autlap$ and
$H/N_2 \cong \GL_d$ canonically.  Moreover, \cite{\GeomIV} exhibits
morphisms
$\rho:  X \rightarrow \grasstd$ and $\sigma: X \rightarrow \toptd$ with
$\phi \circ \rho = R \circ \sigma$ having the following properties: 
$\rho$ is an
$\autlap$-invariant geometric quotient of $X$ by $N_1$, and $\sigma$ a
$\GL_d$-invariant geometric quotient of $X$ by $N_2$.  It is
straightforward to deduce (b). \qed
\enddemo   

The final assertion of Proposition 2.5 yields a twin version of Criterion
2.1 which rephrases existence of a coarse moduli space for the
$d$-dimensional modules with top $T$ in terms of the $\autlap$-space
$\grasstd$.  

\proclaim{Criterion 2.6}  Let $V$ be any $\autlap$-stable subvariety of
$\grasstd$.  Then there exists a coarse moduli space for the
representations $P/C$ with $C \in V$ if and only if the $\autlap$-space
$V$ has a categorical quotient which is an orbit map.  In case of
existence, such a quotient $V / \autlap$ is the moduli space we are
looking for. \qed \endproclaim 

We also glean an alternate characterization of top-stable degenerations
from Proposition 2.5.   Observe that, whenever
$\phi(C) = P/C \cong R(x)$, the correspondence of the theorem pairs the
$\autlap$-orbits contained in the closure of $\aut_\la(P).C$ in
$\grasstd$ with the $\GL_d$-orbits contained in the closure of
$\GL_d.x$ in
$\toptd$.  In other words, the top-stable degenerations of
$P/C$ are precisely the factors $P/C'$ with $C' \in
\overline{\autlap.C}$. 

\subhead 2.D. The partition of $\grasstd$ in terms of radical layers
\endsubhead

In translating geometric features of the variety
$\grasstd$ into representation-theoretic data, the  partition of
$\grasstd$ introduced in this section is pivotal.  Again, we let the top
$T$ equal
$\bigoplus_{1 \le i
\le n} (\la e_i/Je_i)^{t_i}$ and let $P$ be its projective cover; 
moreover, we abbreviate the dimension
$\sum_i t_i$ of $T$ by $t$.   Further, we will systematically identify
isomorphic semisimple modules in the sequel.

\definition{Definitions 2.7} {\bf (i)}  A {\it $d$-dimensional
semisimple sequence with top $T$\/} is a sequence
$$\SS = (\SS_0, \SS_1, \dots, \SS_L)$$ with the following properties:
$\SS_0 = T$ and each $\SS_i$ is a submodule of $J^i P / J^{i+1}P$, such
that $\sum_{i} \dim
\SS_i =d$.  
\smallskip

{\bf (ii)} Given a $d$-dimensional semisimple sequence
$\SS = (\SS_0, \dots, \SS_L)$ with top $T$, we set
$$\grassSS = \{C \in \grasstd \mid \SS(P/C) = \SS \},$$ and call the
$(L+1)$-tuple of integers 
$$\underbardim\SS= (\dim_K \SS_0,
\dots,\dim_K \SS_L)$$ the {\it dimension vector of 
$\SS$\/}.
\enddefinition

Since we are identifying isomorphic semisimple modules, there are only
finitely many $d$-dimensional semisimple sequences. The important
examples are the radical layerings of the
$d$-dimensio\-nal modules with top $T$;  indeed, whenever
$C \in \grasstd$, the sequence $\SS(P/C)$  --  see end of Section 1  --  
meets our criteria.

We start by recording the `layer-stable analogue' of the final
observation of 2.C.  Namely  --  again by Proposition 2.5  --  the
layer-stable degenerations of a module $P/C$ with $C \in \grassSS$ are
precisely the modules represented by the points in the relative closure
of $\autlap.C$ in $\grassSS$.

It is, moreover, obvious that each of the sets $\grassSS$ is stable under
the $\aut_\la(P)$-action and that  $\grasstd$ is the disjoint union of
the
$\grassSS$, where $\SS$ ranges through the $d$-dimensional semisimple
sequences with top $T$.  We point out that, in general, this partition
falls slightly short of being a stratification in the technical sense: 
While the sets $\grassSS$ are locally closed subvarieties of
$\grasstd$   --  see part (ii) of the upcoming lemma  --  their closures
are not always unions of $\grass{\SS'}$'s.  However, this lack is
compensated for by the following asset (see part (i) of the lemma): 
Endow the $(L+1)$-tuples of non-negative integers with the lexicographic
order, and let $\bold{d} \in \NN_{0}^{L+1}$.  Then the union
$\bigcup_{\underbardim \SS >
\bold{d}} \grassSS$, where $\SS$ runs through the 
$d$-dimensional semisimple sequences with top $T$, is closed in
$\grasstd$.  

We can do a little better than that by introducing the following partial
order on the semisimple sequences themselves.  Namely, for two
$d$-dimensional semisimple sequences $\SS = (\SS_0, \dots, \SS_L)$ and
$\SS' = (\SS'_0, \dots, \SS'_L)$, we define
$$\SS \le \SS'$$  if and only if $\bigoplus_{0 \le i \le L} \SS_i =
\bigoplus_{0 \le i \le L} \SS'_i$ and either $\SS = \SS'$ or else
$\SS_j$ is a proper direct summand of
$\SS'_j$ for the smallest index $j$ with $\SS_j \ne \SS'_j$.  Clearly,
this partial order is compatible with the order on the dimension vectors,
in the sense that 
$\SS \le \SS'$ implies $\underbardim \SS \le \underbardim \SS'$.

\proclaim{Lemma 2.8}  

{\rm (i)} Let $\SS$ be any $d$-dimensional semisimple sequence with top
$T$.  Then $\bigcup_{\SS' \le \SS} \grass{\SS'}$ is open in
$\grasstd$, as is the union $\bigcup_{\SS' < \SS} \grass{\SS'}$.
Consequently , the sets $\bigcup_{\underbardim \SS \ge \bold{d}}
\grassSS$ and  $\bigcup_{\underbardim \SS > \bold{d}} \grassSS$ are
closed for any
$(L+1)$-tuple $\bold{d}$ of non-negative integers. 
\smallskip

{\rm (ii)}  Each of the sets $\grassSS$ is locally closed in
$\grasstd$.
\smallskip

{\rm (iii)}  If $M$ is a $d$-dimensional module with top $T$  and
$M'$ a top-stable degeneration of $M$, then
 $\underbardim\SS(M) \le \underbardim\SS(M')$. 
\smallskip
\endproclaim

\demo{Proof}  We fix a sequence $z_1, \dots, z_t$ of top elements of $P$,
where each
$z_r$ is  normed by a primitive idempotent, $e(r)$ say, such that $P=
\bigoplus_{1 \le r \le t} \la z_r$.  

(i) Let $\B$ be a subset of the projective $KQ$-module $\bigoplus_{1 \le
r \le t} KQe(r)$ consisting of elements of the form $p \e(r)$, where $p$
is a path starting in the vertex $e(r)$, and
$\e(r)$ is the element of the direct sum carrying the idempotent $e(r)$
in the
$r$-th slot and $0$ elsewhere.  We say that such a set $\B$ of ``labeled
paths" is {\it admissible\/} if $\B$ includes the candidates 
$\e(1), \dots, \e(r)$  of length
zero and $|\B| = d$.  Moreover, given a $d$-dimensional
semisimple sequence $\SS$ with top $T$, we call an admissible set $\B$ 
{\it compatible with
$\SS$\/} in case, for each $l \in \{0, \dots, L\}$ and each
$i \in \{1, \dots, n\}$, the number of paths of length $l$ in $\B$
ending in $e_i$ equals the multiplicity of the corresponding simple
module $S_i$ in $\SS_l$ (for $l=0$, this is automatic).  As before, $d'
=
\dim P - d$.  For every admissible set $\B$, finally, we let
$\Schu(\B)$ be the set of all points $C$ in the classical
Grassmannian $\Gr(d',JP)$ of $d'$-dimensional subspaces of $JP$ with the
property that 
$$C \oplus \bigl( \sum_{r \le t} \  \sum_{p \e(r) \in \B}  K p z_r
\bigr) = P;$$
here we identify the paths $p$ with their residue classes modulo $I$. 
If
$\Schu(\B)$ is nonempty, then
$\Schu(\B)$ is a Schubert cell of
$\Gr(d',JP)$, and the open subset $\Schu(\B) \cap \grasstd$ of
$\grasstd$ consists precisely of those $d'$-dimensional
$\la$-submodules $C
\subseteq JP$ for which $P/C$ has basis $\bigcup_{1 \le r \le t} \{p
z_r+C \mid p \e(r) \in \B\}$.   

To verify openness of $\bigcup_{\SS' < \SS} \grass{\SS'}$, it now
suffices to notice that 
$\bigcup_{\SS' < \SS} \grass{\SS'}$ is the union of the intersections
$\Schu(\B) \cap \grasstd$, where $\B$ traces the admissible sets that
are compatible with some $\SS'$ smaller than
$\SS$.  The inclusion of the latter union in the former follows from
the fact that, whenever $\B$ is compatible with $\SS'$ and $C \in
\Schu(\B) \cap \grasstd$, the module
$P/C$ has radical layering $\SS(P/C) \le \SS'$; the other inclusion is
obvious.   Analogously, $\bigcup_{\SS' \le \SS}
\grass{\SS'}$ is seen to be open.

(ii)  This follows from (i) and the fact that $\grassSS$ is the
intersection of $\bigcup_{\SS' \le \SS} \grass{\SS'}$ with $\bigcup_{\SS'
\not< \SS} \grass{\SS'}$.  

(iii)  It suffices to observe that the relative closure of
$\grass{\SS(M)}$ in $\grasstd$ is contained in
$\bigcup_{\underbardim\SS\ge \underbardim\SS(M)} \grassSS$.
\qed
\enddemo

\subhead  2.E. General results on the $\autlap$-orbits \endsubhead

The following proposition summarizes the information on the
$\autlap$-orbits of $\grasstd$  --  alias the fibres of the
representation map $\phi$  --  available without any restriction on
$T$.  Strengthened versions for the squarefree case will come to bear in
Section 4. 
  
Again, we let $T = \bigoplus_{1 \le i \le n} (\la e_i/Je_i)^{t_i}$, set 
$t = \sum_i t_i$, and fix a sequence 
$z_1, \dots, z_t$ of top elements of $P$,  each
$z_r$ being normed by a primitive idempotent $e(r)$, such that $P=
\bigoplus_{1 \le r \le t} \la z_r$.  If
$\H$ is the subgroup of $\autlap$ consisting of all automorphisms which
leave the subspace $\bigoplus_{1
\le r \le t} Kz_i$ invariant, and $\U$ the unipotent radical of
$\aut_\la(P)$, then Observation 2.3  tells us that $\U$ is normal in
$\aut_\la(P)$ and
$\aut_\la(P)$ is a semidirect product $\U \rtimes \H$.  The incarnation 
of the maximal torus of $\H$ afforded by our choice of top elements of
$P$ is  
$\T = (K^*)^t$, where we identify $(a_1, \dots, a_t) \in \T$ with the
automorphism of $P$ given by $z_i \mapsto a_i z_i$.  Finally, given
$h \in \H$ and $C\in \grasstd$, we denote by $h.\U.C$ the set $\{hu.C
\mid u \in \U\}$. Due to the normality of $\U$ in $\autlap$, the set
$h.\U.C$ equals the
$\U$-orbit of $h.C$.

\proclaim{Proposition 2.9}  Let $C \in \grasstd$.

{\rm (1)}  The orbit $\autlap.C$ has dimension
$$\dim_K \End_\la(P)  - \dim_K \Hom_\la(P,C) - \dim_K \End_\la(P/C),$$ 
and this dimension is generically constant on the irreducible components
of $\grasstd$. 
\smallskip

{\rm (2)}  The $\U$-orbit $\U.C$ in $\grasstd$ is isomorphic to
$\AA^m$ with 
$$m = m(P/C) = \dim_K \hom_\la(P,JP/C) - \dim_K
\hom_\la(P/C,JP/C).$$   
\smallskip   

{\rm (3)}  There exists a point $C' \in \autlap.C$ such that the
$\T$-orbit $\T.C'$ is isomorphic to the torus
$(K^*)^{t - \s}$, where $\s = \s(P/C)$ is the number of indecomposable
summands of $P/C$.  
\smallskip

{\rm (4)}  The full orbit $\autlap.C$ is the disjoint union of the
subvarieties $h.G.C \cong \AA^{m(P/C)}$, for $h
\in \H$,  and
$\H$ acts as a transitive permutation group on these subvarieties.
\endproclaim

Before proving the proposition, we single out a lemma for repeated
reference.

\proclaim{Lemma 2.10}  For any $C \in \grasstd$, the orbit map
$\aut_\la(P) \rightarrow
\aut_\la(P).C$, $f \mapsto f.C$, is separable, and hence  $\autlap.C$ is
isomorphic to the homogeneous $\autlap$-space $\autlap /
\stab_{\autlap} C$. Moreover, the dimensions of the stabilizers of $C$
in $\autlap$, resp\. in $\U$, are:
 $$\align \dim \stab_{\autlap} C &= \dim_K \Hom_\la(P, C) + \dim_K
\End_\la(P/C) \qquad\quad\text{\it and} \\
\dim \stab_{\U} C &=  \dim_K \Hom_\la(P, C) + \dim_K \Hom_\la(P/C, JP/C).
\endalign$$
\endproclaim

\demo{Proof} We first address the separability claim. Combining
\cite{\Bor, Proposition 6.7} with
\cite{\Kra, AI.5.5, Satz 2}, we see that it suffices to check
reducedness of the schematic fibre of the orbit map, i.e., reducedness
of the stabilizer subgroup $\stab_{\autlap} C$ of $C$.  Reducedness of
this stabilizer, in turn, can be deduced from the fact that it arises as
the solution set of a system of linear equations over $K$; ascertaining
the latter fact  is a matter of routine, if notationally cumbersome.

For the first dimension formula, consider the homomorphism $\rho:
\stab_{\autlap} C
\rightarrow \Aut_\la(P/C)$ of algebraic groups sending any map $f$ in the
stabilizer of $C$ to the induced automorphism of $P/C$.  Clearly, $\rho$
is onto and the kernel is the subgroup of $\autlap$ consisting of the
automorphisms of the form
$\id \, + \ g_0$ with $g_0 \in \Hom_\la(P,C)$.  This proves the first
equality. 

As for the second, consider the restriction of $\rho$ to $\stab_{\U} C$,
and notice that the kernel is the same as before, while $\rho(\stab_{\U}
C)$ equals the unipotent radical of $\Aut_\la(P/C)$, namely the subgroup
consisting of the automorphisms of the form $\id \, + \ g_1$ with
$g_1 \in \Hom_\la(P/C, JP/C)$.  \qed \enddemo   

\demo{Proof of Proposition {\rm 2.9}}  (1) Since
$\dim \autlap.C$ $=$ $\dim \autlap - \dim \stab_{\autlap} C$, the
claimed equality can be derived from Lemma 2.10.  For the generic
behavior of the fibre dimension see, e.g., 
\cite{\Kra, II.2.6}.

(2)  In view of Lemma 2.10, $\U.C$ is isomorphic to the homogeneous space
$\U / \stab_{\U} C$, whence the claim concerning the fibre structure
follows from Rosenlicht's Theorem\cite{\Ros}.  The assertion concerning
the fibre dimension is a consequence of the equality
$\dim \U$ $=$
$\dim_K
\Hom_\la(P,JP)$ and Lemma 2.10. 

(3) Clearly, $\T.C \cong \T/ \stab_{\T} C$ is a torus.  For the
following choice of $C' \in \autlap.C$ we will show
$\dim \stab_{\T} C' = \s$:  Decompose $P/C$ into indecomposable
summands, say $P/C  \cong P_1/C'_1 \oplus \cdots \oplus P_{\s}/C'_{\s}$,
where the $P_i$ are suitable direct summands of $P$;  it is clearly
harmless to assume  $P_i = \bigoplus_{j \in I(i)}
\la z_j$ for some partition $\bigsqcup_{1 \le i \le \s} I(i)$ of $\{1,
\dots, t\}$.  Now set $C'= \bigoplus_{1 \le i \le
\s} C'_i$.   For simplicity of notation, we  assume $C' = C$.  Let
$\widehat{\T}$ be the
$K$-subspace of $\End_\la(P)$ consisting of {\it all\/} endomorphisms of
the form $(a_1, \dots, a_t) \in K^t$, not just the invertible ones;
again we identify $(a_1,
\dots, a_t)$ with the endomorphism of $P$ sending $z_i$ to $a_i z_i$. 
Moreover, $\stab_{\hat{\T}} C$ will stand for the $K$-subspace   of
$\widehat{\T}$ consisting of the endomorphisms 
$f$ with $f(C) \subseteq C$.  Clearly, $\stab_{\T}C$ is a dense open
subvariety of
$\stab_{\hat{\T}} C$, whence the two varieties have the same
dimension.   As for the dimension of the latter, 
$\stab_{\hat{\T}} C$ obviously contains the span of the vectors
$\bold{1}_{I(i)}$, defined as having $j$-th entry $1$ if $j \in I(i)$
and $0$ otherwise.  That these vectors in fact form a basis for
$\stab_{\hat{\T}} C$ is an immediate consequence of the
indecomposability of the $P_i/C_i$.  We infer that $\stab_{\T} C$ has
dimension $\s$, which makes the dimension of
$\T.C$ equal to $t - \s$ as asserted.

(4)  is obvious in view of the remarks preceding the proposition.    \qed
\enddemo

As part (3) of the proposition suggests, the dimension of the stabilizer
subgroup $\stab_{\T} C'$ varies from one point $C' \in \autlap.C$ to
another. In fact, $\stab_{\T} C'$ may be reduced to the scalar multiples
of the identity $\bold{1}$ of $\T$, even when
$\s(P/C)$ is large.

\head 3.  $\grasstd$ for squarefree $T$
\endhead

The diagram below is to remind the reader of the relationships among the
varieties in which we are interested.

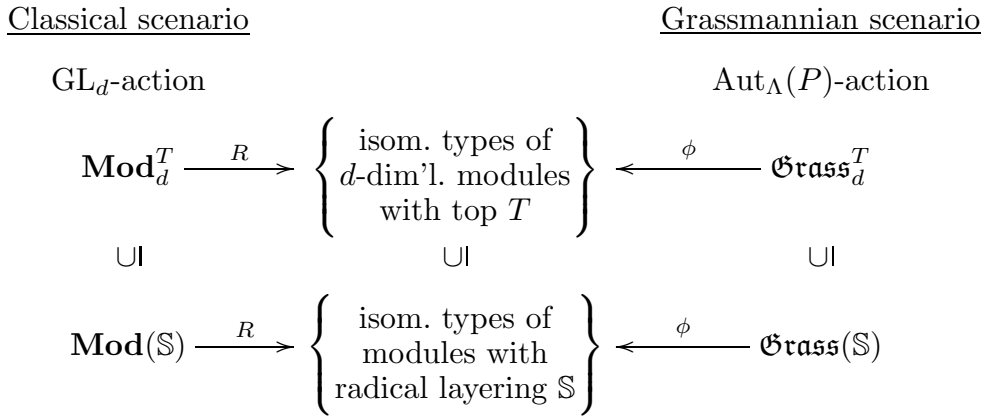
\ignore
$$\xymatrixrowsep{0.3pc}\xymatrixcolsep{0.2pc}
\xymatrix{
\txt{\underbar{Classical scenario}} &&&&&&\txt{\underbar{Grassmannian
scenario}} \\
\txt{$\GL_d$-action} &&&&&&\txt{$\autlap$-action} \\
\toptd \ar[rr]^-{R} &&&\save**\frm{\{}**\frm{\}}\restore
{\txt{isom.~types of\\
$\,d$-dim'l.~modules$\,$\\ with top $T$}}
  &&&\grasstd \ar[ll]_-{\phi} \\
\vsubseteq &&&\vsubseteq &&&\vsubseteq \\
\operatorname{\bold{Mod}}(\SS) \ar[rr]^-{R}
&&&\save**\frm{\{}**\frm{\}}\restore {\txt{isom.~types of\\ modules
with\\ $\,$radical layering $\SS\,$}}
  &&&\grassSS \ar[ll]_-{\phi}
 }$$
\endignore
\medskip

\noindent Throughout this section, $T$ denotes a {\it squarefree\/}
semisimple module, say
$$T = \bigoplus_{1 \le r \le t} \la e_r/Je_r$$  and $P= \bigoplus_{1\le
r\le t} \la e_r$ its projective cover; in particular, $P$ is now a left
ideal of $\la$. We again set $d' = \dim_K
P-d$. Our primary goal is to construct finite
$\autlap$-stable affine open covers of the subvarieties $\grassSS
\subseteq \grasstd$.

\subhead 3.A.  Convenient local coordinates for $\grasstd$ and a
pivotal family of rep\-re\-sen\-ta\-tions para\-me\-trized by $\grasstd$
\endsubhead

We could have postponed the restriction to squarefree $T$ until the next
section, but adopt it here for its notational pay-off.  The coordinatized
open cover of
$\grasstd$ which we introduce first does not consist of $\autlap$-stable
sets in general.  It depends on a fixed choice of top elements of $P$,
whereas the cover at which we are ultimately aiming is $\autlap$-stable
and invariant under such choices.  At the outset, our preferred bases
for factor modules
$P/C$ may appear more specialized than necessary for the present
purpose. They actually are;  yet, it is precisely this special format
which will make them useful in the sequel. 

The second definition below is merely of temporary importance as an
auxiliary to proofs.

\definition{Definitions 3.1} {\bf (i)}  A {\it $d$-dimensional
skeleton with top
$T$\/} is a set $\S$ that consists of $d$ distinct paths of lengths $\le
L$ in $KQ$ and satisfies the following conditions:
\roster

\item"$\bullet$" Each path in $\S$ starts in one of the vertices $e_1,
\dots, e_t$;

\item"$\bullet$" $\S$ is closed under right subpaths; and

\item"$\bullet$" $\S$ contains $\{e_1, \dots, e_t\}$. 
\endroster   
\smallskip

{\bf (ii)}  Given a $d$-dimensional skeleton $\S$ with top
$T$, we define the auxiliary set 
$$\operatorname{Aux}(\S) = \{C \in \grasstd \mid  P/C \ \text{has
basis}\  \{ p  + C \mid p \in \S \}  \}.$$ 
\smallskip  
\enddefinition

 A
$d$-dimensional skeleton
$\S$ with top $T$ can be visualized as an undirected graph with $d$
nodes, namely the disjoint union of $t$ trees,  represented by the sets
$\{p \in \S
\mid p = p e_r \}$ for $1 \le r \le t$.
Clearly, there are only finitely many $d$-dimensional
skeletons.  Moreover, $\operatorname{Aux}(\S)$ is empty  unless the
subset $\S$ of $P$ is linearly independent over
$K$.   On the other hand, each point $C \in \grasstd$ belongs to at
least one of the sets $\operatorname{Aux}(\S)$; in other words, the
$\operatorname{Aux}(\S)$ cover $\grasstd$ --  see Lemma 3.9 for a
sharper statement.

By definition, the sets $\operatorname{Aux}(\S)$ are
intersections of
$\grasstd$ with open Schubert cells in the full Grassmannian
$\operatorname{Gr}(d', JP)$ in which $\grasstd$ is located. 
Namely, for any skeleton
$\S$, the relevant open cell is 
$$\{C \in \operatorname{Gr}(d', JP) \mid C \oplus  \biggl(\sum
\Sb p 
\in
\S \\ \len(p) \ge 1 \endSb K p \biggr)\ \  = \ \ JP\}.$$

\proclaim{Observation 3.2} The $\operatorname{Aux}(\S)$ are open subsets
of
$\grasstd$ \qed
\endproclaim
 
Next, we will use this cover to construct a family of $d$-dimensional
modules with top $T$ parametrized by $\grasstd$.  The reason why this
family is of interest  lies in the fact that it melts down to {\it
the\/} universal family whenever a fine moduli space for our problem
exists.

\proclaim{Lemma 3.3}  There exists a family $(\Gamma, \gamma)$ of
$d$-dimensional modules parametrized by $\grasstd$ such that, for each
$C \in
\grasstd$, the fibre of $\Gamma$ over $C$ (with the $\la$-structure
induced by $\gamma$) is isomorphic to $P/C$. 
\endproclaim

\demo{Proof} Let $\modlad$ be the affine variety of Section 2.A and
$\sigma$ a $d$-dimensional skeleton with top $T$.  We start by
pinning down the obvious morphism
$\rho_\S$ from
$\operatorname{Aux}(\S)$ to $\modlad$ which matches our
coordinatization.  To that end, we assume that $\operatorname{Aux}(\S)
\ne \varnothing$.  This implies that the
$K$-subspace $V$ of
$P$ generated by $\S$ has dimension $d$ and, given any $C \in
\operatorname{Aux}(\S)$, each element of
$P$ has a unique expansion as a
$K$-linear combination of  the elements $p\in\S$ modulo
$C$.  View $V$ as $K^d$, identifying the distinguished basis
$\S$ of $V$ with the canonical basis of
$K^d$, where $\S$ is ordered as follows:  Start by ordering the finite
set $\P$ of all paths of length $\le L$ in $KQ$ in such a fashion that
the paths starting in $e_i$ precede those starting in $e_j$ whenever $i
< j$.  Then let $p_1$ be that element of $\S$ which is minimal in $\P$,
next let $p_2\in\S$ be as small as possible in $\P
\setminus \{p_1\}$, etc.   

Given $C \in \operatorname{Aux}(\S)$, we now define
$$\rho_\S (C)  = \bigl(\rho_\S (C)_{\alpha}\bigr)_{\alpha \in Q^*}$$ 
where $\rho_\S (C)_{\alpha}$ is the $d \times d$ matrix over $K$
recording in its columns the expansion coefficients (relative to $\S$) of
the multiples
$\alpha p$ modulo $C$, for $p \in \S$.   

Thus the trivial bundle $\operatorname{Aux}(\S)
\times K^d$ becomes a family of
$\la$-modules via the algebra homomorphism from $\la$ to the
endomorphism ring of
$\operatorname{Aux}(\S) \times K^d$, defined by 
$$\alpha \mapsto \bigl( (C,v) \mapsto (C, \rho_\S (C)_{\alpha} v)
\bigr).$$   The appropriate choice of morphisms 
$g_{\sigma, \sigma'}: \operatorname{Aux}(\S) \cap \operatorname{Aux}(\S')
\rightarrow \GL_d$ to glue the $\operatorname{Aux}(\S) \times K^d$
together to a vector bundle over $\grasstd$ is now automatic:  Namely
just send each $C \in \operatorname{Aux}(\S)
\cap  \operatorname{Aux}(\S')$ to the transition matrix recording the
coefficients of the $p' + C$ (for  
$p'\in \sigma'$) relative to the $p + C$, in the prescribed order; then,
clearly, the $g_{\sigma, \sigma'}$ satisfy the relevant cocycle
condition.   The resulting vector bundle over $\grasstd$ will be denoted
by
$\Gamma$.  Since our algebra homomorphisms $\la
\rightarrow \End(\operatorname{Aux}(\S)\times K^d)$ are compatible with
this gluing, they yield an algebra homomorphism $\gamma: \la
\rightarrow \End(\Gamma)$ and thus a family of $d$-dimensional
$\la$-modules as required. \qed \enddemo

\subhead  3.B.  The distinguished $\autlap$-stable affine open cover of
$\grassSS$
\endsubhead

The key to the announced cover comes from the concrete counterparts of
the skeletons of section 3.A.  By construction, this new cover
will depend only on the quiver of $\la$ and the ideal $I$ of relations. 
In fact, up to permutation and birational equivalence of the irreducible
components, it is an isomorphism invariant of $\la$ (see 3.18).

\definition{Definitions 3.4} {\bf (i)} Let $M$ be a $d$-dimensional
module with top $T$.
 A $d$-dimensional skeleton $\S$ (with top $T$) is called a
{\it skeleton of
$M$\/} provided there exists a sequence $x_1, \dots, x_t$ of top
elements of $M$ with $x_r = e_r x_r$ such that, for every integer $l$
between $0$ and $L$, the union  
$$\bigcup_{1 \le r \le t} \{ p x_r +J^{l+1}M \mid p = p e_r \in
\S,\
\len(p) = l\}$$    is a $K$-basis for $J^lM/J^{l+1}M$. 
\smallskip

{\bf (ii)}  For any $d$-dimensional skeleton $\S$ with top $T$,
we set
$$\grassS = \{ C \in \grasstd \mid \S \ \text{is a skeleton of\ } P/C
\}.$$
\smallskip

{\bf (iii)}  Given a skeleton $\S$ and a semisimple sequence
$\SS$, we say that $\S$ is {\it compatible with $\SS$\/} in case, for
each $i\in \{1,\dots,n\}$, the number of paths of any prescribed length
$l$ in $\S$ which end in the vertex $e_i$ equals $\dim_K e_i \SS_l$.
\enddefinition

First we observe that the notion of a skeleton of a module $M$ as in (i)
does not depend on a choice of top elements of $M$.  Indeed, as is
readily confirmed, if $\S$ is a skeleton of
$M$, then any choice of top elements $y_1 = e_1 y_1, \dots, y_t = e_t
y_t$  of $M$ yields $K$-bases 
$\bigcup_{1 \le r \le t} \{ p y_r+J^{l+1}M
\mid p = p e_r \in \S,\, \len(p) = l \}$ for the radical layers
$J^lM/J^{l+1}M$. 

It is obvious that isomorphic modules have identical sets of skeletons,
which means that the sets
$\grassS$ are stable under the action of $\autlap$.  Moreover, we notice
that, whenever $C$ belongs to $\grassS$, the skeleton $\S$ is compatible
with the sequence
$\SS(P/C)$ of radical layers of $P/C$.  Hence $\grassS \cap \grassSS =
\varnothing$ if $\S$ fails to be compatible with $\SS$.  On the other
hand, $\grassS$ is contained in $\grassSS$, whenever $\S$ is compatible
with $\SS$.  Further properties of the sets  $\grassS$ are collected in
the following theorem.  In particular, we glean from it that the
$\grassS$ are locally closed affine subvarieties of $\grasstd$.  
 
\proclaim{Theorem 3.5}  Let $\SS$ be a $d$-dimensional semisimple
sequence with top $T$. 

The sets $\grassS$, where $\S$ runs through the skeletons
compatible with $\SS$, form a finite $\autlap$-stable affine open cover
of the quasi-projective variety
$\grassSS$.  A family of polynomials determining any of the
$\grassS$ in a suitable affine coordinate system can be
{\rm{(}}algorithmically{\rm{)}} obtained from the quiver
$Q$ and a finite set of relations generating the ideal $I$ as a left
ideal of $KQ$.   \endproclaim  

The fact that the sets $\grassS$ form an open cover of $\grassSS$ will
be proved in Lemmas 3.8 and 3.9; we emphasize that the
$\grassS$ fail to be open in $\grasstd$ though. The remaining claim of
Theorem 3.5,  namely that each of the sets
$\grassS$ is an affine variety that can be obtained from quiver and
relations of $\la$, is less obvious.  We will tackle it in several
steps, providing further detail along the way.   
 
\definition{ Remark 3.6}    The affine cover of Theorem 3.5 permits to
resolve the isomorphism problem for
$d$-dimensional top-$T$ modules  with low computational investment: 
Given any point $C$ in
$\grassS$, expressed in terms of an affine coordinate system for
$\grassS$, the
$\aut_\la(P)$-orbit of $C$ can be obtained from a system of at most
$d \cdot (\dim_K J/J^2)$ linear equations (see \cite{\preview}).
\enddefinition

The following elementary example shows that, if $T$ contains squares of
simples, the variety $\grassSS$ will not have an $\autlap$-stable affine
cover in general.  In fact, we exhibit an $\autlap$-orbit that fails to
be quasi-affine for a situation where $T$ is a square.

\definition{Example 3.7}  Let $\la$ be the representation-finite
hereditary algebra $KQ$, where $Q$ is the quiver

\ignore
$$\xymatrixrowsep{2pc}\xymatrixcolsep{3pc}
\xymatrix{ 2 &1 \ar[l]_{\alpha} \ar[r]^{\beta} &3
 }$$
\endignore 

If $\SS = (S_1^2, S_2 \oplus S_3)$, then $\grassSS$ is isomorphic to
$\PP^1 \times \PP^1$ and consists of two $\autlap$-orbits, namely the
orbits representing the modules $\la e_1/ \la \alpha e_1
\oplus \la e_1/ \la \beta e_1$ and $\la e_1 \oplus S_1$.  The orbit of
the latter is isomorphic to $\PP^1$. 
\qed\enddefinition    

We return to our blanket hypothesis that $T$ be squarefree.  In
justifying Theorem 3.5, we start with openness of the
$\grassS$ in the pertinent $\grassSS$.

\proclaim{Lemma 3.8}  For any skeleton $\S$ compatible with $\SS$, the
subset
$\grassS$ is open in $\grassSS$. \endproclaim

\demo{Proof}  From Definitions 3.1 and 3.4 one readily deduces that,
under our compatibility hypothesis,
$$\grassS = \operatorname{Aux}(\S) \cap \grassSS.$$  
Observation 3.2 thus proves our claim.
\qed
\enddemo  

 That the $\grassS$ cover $\grassSS$, is guaranteed by the following
elementary lemma. 

\proclaim{ Lemma 3.9}  Every module $M$ with top $T$ has at least one
skeleton. 
\endproclaim

\demo{Proof} First, we set $\S_0 = \{e_1, \dots, e_t\}$  and choose top
elements $x_1=e_1x_1,\dots,x_t=e_tx_t$ of $M$ giving rise to a basis of
$M/JM$.  Next, we pick a set
$\S_1$ of arrows
$p_1$, each starting in one of the vertices in $\S_0$, such that the set
$\{p_1 x_r \mid  p_1\in \S_1\cap KQ e_r,\ 1\le r\le t \}$ yields a
$K$-basis for
$JM$ modulo $J^2M$. The following stage is to choose a set $\S_2$ of
paths
$p_2$ of length 2, each containing some arrow in $\S_1$ as a right
subpath, such that the set
$\{ p_2 x_r \mid p_2\in \S_2 \cap KQ e_r,\ 1\le r\le t\}$ gives
rise to a $K$-basis of
$J^2M/J^3M$.  We continue in this fashion and set $\S= \bigcup_{0\le  l
\le L} \S_l$. (Keep in mind that $J^{L+1}=0$ by our choice of $L$.) \qed
\enddemo

From the proof of Lemma 3.9, we can actually glean the following
stronger statement: Given any path
$p$ of length
$l$  with $p M
\not\subseteq J^{l+1}M$, there exists a skeleton of $M$ which contains
$p$.

We next introduce a convenient affine coordinate grid for each
$\grassS$, and then follow with a family of polynomials derived from the
relations of
$\la$.  That the vanishing set of these polynomials coincides with
$\grassS$, up to an isomorphism which `respects' the   representation
map, is spelled out under the heading of the more precise Theorem 3.14
below.

\definition{3.10. The coordinate system for $\grassS$} Let $\S$ be
a $d$-dimensional skeleton with top $T$.  As we know from 3.A, $\grassS$
embeds into the intersection 
$\operatorname{Aux}(\S)$ of $\grasstd$ with
the open Schubert cell consisting of those subspaces of $JP$ which have
zero intersection with the span of $\S$ (see Definitions 3.1).  In
essence, the coordinates with which we will work are Pl\"ucker
coordinates relative to an ordered basis of
$JP$ which supplements the paths of positive length in
$\S$ to a
$K$-basis of $JP$ made up of paths.  However, our situation permits
restriction to a small subfamily of these coordinates, determined by
paths of the form
$\alpha p$, where
$(\alpha,p)$ is a $\S$-critical pair as defined below.  As
pointed out in the introduction, this economy is crucial in the
analysis of concrete instances; polynomial equations governing these
pared-down coordinates can be read off from quiver and relations of
$\la$.  

A pair $(\alpha,p)$, consisting
of an arrow $\alpha$ and a path $p\in\S$, is called $\S$-{\it
critical\/} in case $\alpha p$ is a path in $KQ$ which does {\it not\/}
belong to
$\S$. 
\enddefinition

Clearly, the isomorphism type of any module $M = P/C$ with skeleton $\S$
is completely pinned down by the (unique) scalars $c_{\alpha p,q}$
appearing in the following equations in $M$:
$$\alpha p + C = \sum_{q \in \S} c_{\alpha p, q} q  + C$$ for the
$\S$-critical pairs $(\alpha, p)$. Moreover, to obtain a more economical
coordinatization, we observe that the coefficient
$c_{\alpha p, q}$ in the above equation is zero unless $q$ is strictly
longer than $p$, and has the same starting and ending vertices as
$\alpha p$.  This motivates the following notation:  For any
$\S$-critical pair
$(\alpha, p)$, we let $\S(\alpha, p)$ be the set of all paths in $\S$
which are at least as long as $\alpha p$ and have the same
end vertices as $\alpha p$.

As a consequence, in pinning down $M$, we only need to keep track of the
scalars arising in the equations 
$$\alpha p + C = \sum_{q \in \S(\alpha,p)} c_{\alpha p, q} q  + C$$ for
the $\S$-critical pairs $(\alpha, p)$.   In other words, as a
$\la$-submodule of $P$, any $C \in \grassS$ is generated by the
differences 
$$\alpha p \ - \sum_{q \in \S(\alpha, p)} c_{\alpha p, q} q.$$  As long
as we keep the skeleton
$\S$ fixed, it is therefore justified to identify the points
$C \in \grassS$ with the corresponding families  
$$\bigl(c_{\alpha p, q} \bigr) \Sb (\alpha, p)\ \S\text{-critical} \\ q
\in \S(\alpha, p) \endSb$$ of scalars in the appropriate affine space
$\AA^N$, where $N$ is the disjoint union of the sets
$\S(\alpha,p)$.       

\definition{3.11. The congruence relation induced by $\S$} Keeping $\S$
fixed, we consider the polynomial ring  
$$\A = \A(\S) := KQ[X_{\alpha p, q} \mid (\alpha, p)\
\S\text{-critical},\ q \in
\S(\alpha,p)]$$ over the path algebra $KQ$.
On the ring $\A$, we consider congruence modulo the left ideal
$$\C = \C(\sigma) := \bigoplus_{t+1\le i\le n} \A e_i \ \ + \
\sum \Sb (\alpha,p)\ \S\text{-critical} \endSb \A \bigl( \alpha p \ -
\sum \Sb q\in \S(\alpha,p) \endSb  X_{\alpha p, q} q
\bigr)$$  and denote this congruence relation by $\seq$.
\enddefinition

The argument backing the following proposition is constructive.

\proclaim{Proposition 3.12}  The quotient $\A/\C$ is a free left module
over  the commutative polynomial ring
$K[ X_{\alpha p, q}]$, with the cosets of the paths in $\S$ forming
a basis.

In other words:  For every $z \in \A$, there is a unique family
of polynomials
$\tau_q(X) = \tau^{z}_q(X)$ in 
$K[ X_{\alpha p, q}]$ such that 
$$z \seq \sum_{q \in \S}
\tau_q(X)\ q.$$ 
\endproclaim

\demo{Proof}  We only verify existence of the $\tau_q(X)$ and leave
uniqueness to the reader.  It is clearly innocuous to assume that
$z$ is a path in
$KQ$ starting in one of the vertices $e_1,\dots,e_t$.  Let $p$ be the
longest right subpath of
$z$ belonging to
$\S$.  We will prove our claim by induction on $\len(z) - \len(p)$.  If
this difference is zero, we are done.  So suppose it is positive and let
$\alpha$ be the (unique) arrow with the property that $\alpha p$ is
again a right subpath of $z$, say $z = z' \alpha p$ for a suitable path
$z'$.  Then $(\alpha,p)$ is a $\S$-critical pair, and
$$z \seq z'\bigl(\sum_{q \in \S(\alpha,p)} X_{\alpha p,q} q\ \bigr) = 
\sum_{q \in \S(\alpha,p)} X_{\alpha p,q} z' q.$$   Since all of the
differences $\len(z'q) - \len(q)$ arising from this equation are
strictly smaller than $\len(z) -
\len(p)$, our induction hypothesis guarantees that each of the summands
$X_{\alpha p,q} z'q$ has the desired form up to $\seq$. Consequently, so
does $z$. \qquad \qed \enddemo  

\definition{3.13. Polynomials for $\grassS$}  Let $\I(\sigma)$ be the
ideal of the polynomial ring $\A(\sigma)$ generated by the ideal $I$
of relations of $\la$, i\.e\., $I(\sigma) = I[X_{\alpha p,q}]$. Moreover,
let
$\K(\sigma)$ be the kernel of the $K[ X_{\alpha p, q}]$-resolution
$$\bigl(K[ X_{\alpha p, q}]\bigr)^\sigma \rightarrow \A(\sigma)/
\bigl(\I(\sigma) + \C(\sigma) \bigr),$$
sending $\bigl(\tau_q(X) \bigr)_{q \in \sigma}$ to $\sum_{q \in \sigma}
\tau_q(X)\overline{q}$; here the cyclic left $\A(\sigma)$-module 
$\A(\sigma)/ \bigl(\I(\sigma) + \C(\sigma) \bigr)$ carries
the inherited $K[ X_{\alpha p, q}]$-structure. Then, as Theorem 3.14
will show, the variety
 $\grassS$ is isomorphic to the vanishing set in $\AA^N =
\AA^{N(\sigma)}$ of the ideal
$\sum_{q \in \sigma} \pi_q\bigl(\K(\sigma) \bigr)$ in the commutative
polynomial ring $K[ X_{\alpha p, q}]$, where the
$\pi_q$ for $q \in \sigma$ are the canonical projections of $(K[
X_{\alpha p, q}]\bigr)^\sigma$. 

Proposition 3.12 allows us to rephrase this in the following leaner
terms: Let
$\R$ be any finite generating set for the left ideal $Ie_1 + \cdots
+Ie_t$ of
$KQ$  --  note that such a generating set always exists since all paths
of lengths $\ge L+1$ belong to $I$.  For each
$\rho\in \R$, let
$\tau^{\rho}_q(X)$ be the unique polynomials in $K[ X_{\alpha p, q}]$
with
$$\rho \  \seq \ \sum_{q\in \S} \tau^{\rho}_q(X)\ q$$ as guaranteed by
Proposition 3.12.  Then the locus
$V\bigl(\tau^{\rho}_q(X) \mid \rho \in \R,\, q \in \sigma \bigr)$ of
these polynomials is an alternate incarnation of $\grassS$, as we will
see.

\enddefinition

 Solely for the purpose of establishing a suitable isomorphism of
varieties, we let $\AffS$  stand for the subvariety of $\AA^N$ defined
by the polynomials
$\tau^{\rho}_q(X)$ for the moment. 

It is straightforward that $\AffS$ is independent of our choice of
$\R$. (In fact, so is the ideal of $K[X_{\alpha p,q}]$ generated by the
$\tau_q^\rho(X)$.)  Furthermore, we emphasize that, in order to proceed,
we need not preselect `eligible' 
$d$-dimensional skeletons, i.e., skeletons which actually occur
as skeletons of
$\la$-modules; 
$\Aff(\S)$ is automatically empty if $\grassS$ is, as is for instance
the case when, for some integer $l$, the set of paths of length $l$ in
$\S$ is linearly dependent modulo $J^{l+1} P$.

\proclaim{Theorem 3.14} Suppose $\S$ is a $d$-dimensional skeleton. There exists an isomorphism  $\psi_\S : \AffS
\rightarrow \grassS$ of varieties, together with a surjective map
$\chi_\S$, such that the following diagram commutes:

\ignore
$$\xymatrixrowsep{2pc}\xymatrixcolsep{4pc}
\xymatrix{
\AffS \ar[dd]^{\cong}_{\psi_{\S}} \ar[dr]^{\chi_{\S}}\\
 &\save+<25ex,0ex> \drop{\txt{ \{\rm isom.~types of
$d$-dim'l.~modules with skeleton $\S$\} }} \restore\\
\grassS \ar[ur]_{\phi_{\S}} }$$
\endignore

\noindent Here $\phi_\S$ is the restriction of the representation map
$\phi$ to $\grassS$, and the maps $\psi_\S$ and $\chi_\S$ are as follows:
$\psi_\S$ sends any point $c = \bigl( c_{\alpha p,q} \bigr) \Sb
(\alpha,p)\
\S\text{-critical} \\ q\in \S(\alpha,p) \endSb$ from $\AffS$ to the
submodule
$$U(c)\ :=\ 
\sum \Sb (\alpha,p)\ \S\text{-critical} \endSb \la \bigl( \alpha p \  -
\sum \Sb q\in \S(\alpha,p) \endSb  c_{\alpha p, q} q
\bigr)$$ of $P = \bigoplus_{1 \le r \le t} \la e_r$, and $\chi_\S$ sends
$c$ to the isomorphism type of the factor module $P/U(c)$.
\endproclaim

We defer the somewhat technical proof of Theorem 3.14 to Section 6;
clearly, Theorem 3.14 covers the remaining claims of Theorem 3.5.

In the sequel, we will identify $\AffS$ with $\grassS$ and $\chi_\S$
with $\phi_\S$.  This is warranted in light of the commutative diagram
linking the pertinent representation maps from
$\grassS$ and $\AffS$ to the modules with skeleton $\S$.  We remark that
$\psi_\S^{-1}$ is a closed immersion of $\grassS$ into the affine space
$\AA^N$.        

An algorithm for finding  polynomials that determine $\grassS$ within
the affine space $\AA^N$ is implicit in the preceding discussion, but
will not be put on a systematic basis here.  We will only illustrate
the procedure with an example.  The following observation
significantly reduces the computational effort  --  we include it, since
we will require the underlying concept of a `route' in the
proof of Theorem 3.14.

\definition{3.15. A crucial property of the left ideal $\C(\sigma)$ of
3.11} Suppose
$u$ is a path of length $l$ which passes through the vertices $(e(0),
e(1),\allowmathbreak \dots,\allowmathbreak e(l))$ in that order; in
particular, this means that
$u$ starts in $e(0)$ and ends in $e(l)$.  Given a
$d$-dimensional skeleton $\S$ with top $T$, such a path
$u$ will be called a {\it route on $\S$\/} if there exist paths $p_0,
p_1, \dots, p_l$ in $\S$ with $\len(p_0) = 0$ and
$\len(p_{i-1}) < \len(p_{i})$ such that each $p_i$ ends in $e(i)$.  Note
that any path in $\S$ is a route on
$\S$, since skeletons are closed under right subpaths by definition.

If $u$ is a path which fails to be a route on $\sigma$, then $u \in
\C(\sigma)$; in other words, $u \seq 0$ under the congruence
relation of 3.11.  Consequently, $u$ can be ignored in the substitution
process leading to polynomials defining $\grassS$.

 To justify this claim, let $C$ be any point in $\grassS$, and suppose
$u = \alpha_s \cdots \alpha_1$ is a path, concatenated from arrows 
$\alpha_i$, such that $u + C$ is nonzero in $M = P/C$.  If, for $j
\le s$, the integer $\l_j$  is such that
$\alpha_j
\cdots \alpha_1  + C \in J^{l_j}M \setminus J^{l_j + 1}M$, then
$1 < l_1 < l_2 < \dots < l_s$; moreover, for each $j$, the set $\S$
contains a path of length $l_j$ ending in the same vertex as
$\alpha_j$.  This shows the path $u = \alpha_s \cdots \alpha_1$ to be a
route on $\S$, thus backing our remark. 
\enddefinition

\definition{Example 3.16}  Consider $\la = KQ /I$, where $Q$ is the
quiver 

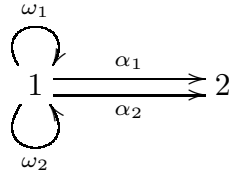
\ignore
$$\xymatrixrowsep{2pc}\xymatrixcolsep{4pc}
\xymatrix{ 1 \ar@(ul,ur)^{\omega_1} \ar@(dl,dr)_{\omega_2}
\ar[r]<0.6ex>^{\alpha_1} \ar[r]<-0.6ex>_{\alpha_2} &2 }$$
\endignore

\noindent and $I \subseteq KQ$ the ideal generated by the four paths
$\omega_i \omega_j$ for $i,j \in \{1,2\}$ and $\alpha_1 \omega_1 -
\alpha_2 \omega_2$.  Clearly, the listed relations also generate $Ie_1$
as a left ideal of $KQ$.

For $T = S_1$ and $d = 4$, we consider the skeleton $\S =
\{e_1,
\omega_1,
\alpha_1
\omega_1, \alpha_2\}$, and determine $\grassS$ as follows. 
First we list the
$\S$-critical pairs $(\alpha,p)$, together with their sets
$\S(\alpha,p)$.  Obviously, there is no harm in omitting pairs
$(\alpha, p)$ with $\alpha p \in I$; the pairs that are left are
$(\omega_2, e_1)$ with $\S(\omega_2, e_1) = \{\omega_1\}$, 
$(\alpha_1, e_1)$ with $\S(\alpha_1, e_1)  = \{\alpha_2, \alpha_1
\omega_1\}$,  and $(\alpha_2, \omega_1)$ with
$\S(\alpha_2, \omega_1) = \{\alpha_1 \omega_1\}$.  This leads to the
following congruences in $KQ[X_1, \dots, X_4]$, where
$X_1,\dots, X_4$ stand for $X_{\omega_2,\omega_1}$,
$X_{\alpha_1,\alpha_2}$, $X_{\alpha_1,\alpha_1\omega_1}$,  $X_{\alpha_2
\omega_1, \alpha_1 \omega_1}$.  Namely, 
$\omega_2 \seq X_1 \omega_1$, $\alpha_1 \seq X_2 \alpha_2 + X_3
\alpha_1 \omega_1$, and
$\alpha_2 \omega_1 \seq X_4 \alpha_1
\omega_1$.   Since all of the four paths $\omega_i \omega_j$ in $I$ fail
to be routes on $\S$, Remark 3.15 tells us that $\omega_i \omega_j
\seq 0$ is automatic, whence the relations $\omega_i \omega_j = 0$ in
$\la$ do not produce any conditions on the $X_k$.  On the other hand,
inserting the above `basic congruences' into the  final relation yields
$\alpha_1 \omega_1 - \alpha_2 \omega_2$ $\seq$ 
$\alpha_1 \omega_1 - X_1 \alpha_2 \omega_1$ $\seq$ $\alpha_1 \omega_1 - 
X_1 X_4 \alpha_1 \omega_1$ $=$ $(1 - X_1X_4)
\alpha_1 \omega_1$.  The last congruence is the expansion of $\alpha_1
\omega_1 - \alpha_2 \omega_2 \in A$ as in Proposition 3.12, whence 
$$\grassS = \{(c_i)_{i \le 4} \in \AA^4 \mid 1 - c_1c_4 = 0\}.
\qquad\square$$    
\enddefinition

\definition{ 3.17.  Functoriality of the $\autlap$-stable affine
charts}  Fix the quiver
$Q$, and identify $T$ with the corresponding set of vertices $\{e_1,
\dots, e_t\}$. Moreover, let
$\I$ be the category of all admissible ideals of the path algebra $KQ$ --
the morphisms in this category are the inclusion maps  --  and define
$\M$ to be the category of left modules over the ring
$\prod_{\sigma} \A(\sigma)$, where
$\sigma$ traces the skeletons with top $T$.  Now define
$F$ to be the functor $\I \rightarrow \M$, which sends any ideal
$I$ in $\I$ to the left module $\prod_\sigma \A(\sigma)/\bigl(I(\sigma)
+ \C(\sigma) \bigr)$ (c\.f\. 3.11 and 3.13).  By the
preceding theorem, the functor $F$ then carries full information about
the representation theory with top $\{1.
\dots, t\}$ of arbitrary basic $K$-algebras with quiver $Q$.  Indeed,
according to 3.13 and 3.14, each image under
$F$ determines a family $\bigl(\grass{I, \sigma}\bigr)_\sigma$, where
$\grass{I,
\sigma}$ denotes the affine chart in 
$\AA^{N(\sigma)}$ which represents the
$KQ/I$-modules with top $T$ and skeleton $\sigma$.  Finally, if $\Var^T$
stands for the category consisting of all families of affine varieties
over $K$ indexed by the skeletons $\sigma$ with top $T$  --  the maps in
this latter category are the corresponding families of morphisms  -- 
then the assignment $\bigl(I(\sigma)
\bigr)_\sigma
\mapsto  \bigl(\grass{I, \sigma} \bigr)_\sigma$ defines a
contravariant functor from $\I$ to $\Var^T$ which reverses inclusions.   
\enddefinition   
       
\definition{3.18.  Uniqueness of the $\autlap$-stable affine charts}
Clearly, the locally closed subsets
$\grassSS$ of
$\grasstd$ are uniquely determined, up to isomorphism, by the
isomorphism type of
$\la$.  The affine open subsets $\grassS$ of any given $\grassSS$, in
their turn, enjoy a somewhat weaker uniqueness property.  Namely, they
are invariant under algebra isomorphisms of $\la$, up to birational
equivalence.  More precisely: Suppose $\la \cong \la'$, where
$\la' = KQ'/ I'$ (it is well-known that
$Q$ and
$Q'$ are isomorphic as directed graphs in this situation). Let
$\SS$ be a semisimple sequence over $\la$, and
$\SS'$ the semisimple sequence over $\la'$ corresponding to $\SS$ under
some algebra isomorphism from $\la$ to $\la'$. Moreover, let
$\Irr(\SS)$ be the set of all irreducible components of the affine
varieties $\grassS$, where
$\S$ runs through the skeletons compatible with $\SS$, and define
$\Irr(\SS')$ analogously.  Then there is a bijection $\Irr(\SS)
\rightarrow \Irr(\SS')$ such that the partners under this pairing   are
birationally equivalent.  One proves this as in the scenario of
varieties of uniserial modules in
\cite{\GeomII}. Finally, one observes that `birationally equivalent' can
be strengthened to `isomorphic' in our uniqueness statement, provided
that $Q$ is without double arrows; indeed, in that case, both the quiver
and the corresponding relations are essentially pinned down by the
isomorphism type of $\la$.  
\enddefinition

\head 4. The moduli problem \endhead

We continue to assume that $T = \bigoplus_{1 \le r \le t} \la e_r /
Je_r$ is squarefree, and $P= \bigoplus_{1\le r\le t}
\la e_r$. The first subsection contains our main results, while in 4.B
we briefly discuss the moduli problem for classes of representations
with fixed radical layering.

\subhead 4.A. Top-stable degenerations and quotients of $\grasstd$
modulo $\autlap$ \endsubhead

From Proposition 2.9 and Theorem 3.5, we derive the following 
consequences concerning top-stable and layer-stable degenerations of
modules with squarefree tops. 

\proclaim{Proposition 4.1}  

{\rm (1)}  All $\autlap$-orbits are quasi-affine subvarieties of
$\grasstd$.
\smallskip

{\rm (2)} Whenever $M$ is a module with top $T$, say $M = P/C$ with $C
\in \grassSS$, the orbit
$\U.C$ of
$C$ under the action of the unipotent radical $\U$ of $\autlap$ is
closed in
$\grassSS$.  It is an affine space $\AA^m$ of dimension 
$$\align m= m(M) &= \mu_T(M) - \bigl(t + \dim_K \Hom_\la(M, JM) \bigr)\\
&= \dim_K \hom_\la(P,JM) - \dim_K \hom_\la(M,JM),\endalign$$ where
$\mu_T(M)$ is the sum of the multiplicities of the simple summands of
$T$ as composition factors of $M$. 
\smallskip 

{\rm (3)} In case  $T$ is simple and $M$ as in part {\rm (2)}, the orbit
$\autlap.C$ is closed in $\grassSS$ and isomorphic to  
$\AA^m$ with 
$$m = \mu_T(M) - \dim_K \End_\la(M),$$ where $\mu_\T(M)$ is the
multiplicity of $T$ in $M$.  
\endproclaim

\noindent {\it Remark.\/}  In \cite{\preview}, it will be shown that,
for squarefree $T$, all $\autlap$-orbits of $\grasstd$ are direct
products of affine spaces and tori as follows:  $\autlap.C \cong
\AA^m \times (K^*)^r$, where $m$ is as under (2) and $r$ equals
the difference $\dim T$ minus the number of indecomposable summands of
$M$. 

\demo{Proof of Proposition 4.1}  The first assertion is an immediate
consequence of the existence of an $\autlap$-stable affine cover of
$\grasstd$ (Theorem 3.5).  As for the second, let $C \in \grassSS$.  To
see that the $\U$-orbit of $C$ is closed in $\grassSS$, we
need only make sure that its intersections with the open subvarieties
$\grassS$, where $\S$ runs through the skeletons compatible with $\SS$,
are relatively closed.  But in light of Theorem 3.5, this follows from a
result of Kostant and Rosenlicht which guarantees that the orbits of a
unipotent group morphically acting on an affine variety are closed (cf\.
\cite{\RosII, Theorem 2}, and \cite{\Hum, p.115, Exercise 8} for a
simple alternate proof).  The claimed dimension of
$\U.C$ can be deduced from Proposition 2.9. For part (3)
finally, observe that $\autlap \cong K^* \times
\U$, with $K^*$ acting trivially, and apply part (2).  
 \qed
\enddemo  

A module is called {\it local\/} in case it has simple top;
thus the local modules are characterized by their having
precisely one maximal submodule.

\proclaim{Theorem 4.2}  Suppose that $M$ is a $d$-dimensional module
with top $T$, say $M \cong P/C$ with $C \in
\grasstd$.  Then the following conditions are equivalent:
\roster
\item"(1)" $M$ has no proper top-stable degenerations.

\item"(2)"  The orbit $\autlap.C$ is closed in $\grasstd$.

\item"(3)"  The orbit $\autlap.C$ is a singleton.

\item"(4)"  $C \subseteq P$ is invariant under
$\la$-endo\-morph\-isms of $P$, i\.e\., $C u \subseteq C$ for any path
$u$ starting in one of the vertices $e_r$ with $r \le t$.   
 
\item"(5)" $M$ is a direct sum of local modules and $$\mu_\T(M) = t +
\dim_K
\Hom_\la(M, JM),$$ where $\mu_T(M)$ is the sum of the multiplicites of
the simple summands of $T$ as composition factors of $M$.    
\endroster
 
\endproclaim  

\demo{Proof} The equivalence of (1) and (2) is immediate
from Proposition 2.5, and the implications
(4)$\implies$(3)$\implies$(2) are trivial.

To verify (2)$\implies$(3), suppose (2) holds.   Using the fact that
$\grasstd$ is a projective variety, we infer completeness of
$\autlap.C$.  On the other hand, this orbit is quasi-affine by
Proposition 4.1, and therefore must be zero dimensional.  Since
$\autlap.C$ is irreducible, (3) follows. 

For (3)$\implies$(4), assume (3), which means that $C$ is invariant under
automorphisms of $P$.  Given any $f \in \End_\la(P)$, pick a nonzero
scalar $a$ such that $-a$ is not an eigenvalue of $f$.  Then $a\cdot\id
+f$ belongs to $\autlap$, and hence leaves $C$ invariant; but this means
$f(C) \subseteq C$.

(3)$\implies$(5).  By Proposition 2.9(3), the number of summands in a
decomposition of $M$ into indecomposables equals $t$.  Since $t = \dim_K
T$, this means that $M$ splits into local summands.  The numerical
equality follows from Proposition 4.1(2).

(5)$\implies$(3).  Without loss of generality, we may assume that $C$
equals the point $C' \in \autlap.C$ specified in Proposition 2.9(3). 
Thus $\T.C = \{C\}$.  Moreover, (5) guarantees that
$\U.C \cong \AA^0$, where $\U$ again denotes the unipotent radical
of $\autlap$.  Hence 
$\autlap.C$ is a singleton.
\qed
\enddemo

In case $T$ contains squares, the statements of the theorem fail in
general.  For instance, in Example 3.7, let
$M = \la e_1 \oplus S_1$.  Then $M= P/C$ with $C= \la(0,\alpha)+
\la(0,\beta)$ has no proper top-stable degeneration, while $\autlap.C$ 
is the projective line.

The simplest instance of a (local) module $M$ having a proper top-stable
degeneration can be found in  Example 2.4.  In this example,
$\grasstd$ consists of two $\autlap$-orbits, those of the modules $M =
\la e_1/ \la
\alpha$ and $M' = \la e_1 / \la \alpha \omega$.  The $\autlap$-orbit of
$C =   \la
\alpha$ is isomorphic to $\AA^1$, while that of $C' =  \la \alpha
\omega$ is a singleton.  Thus $M$ degenerates top-stably to $M'$. 
However, $M'$ is not a layer-stable degeneration of $M$, as is predicted
by the following consequence of Proposition 4.1(3).   

\proclaim{Corollary 4.3}  If $M$ is a $\la$-module with simple top, then
$M$ has no proper layer-stable degenerations. \qed \endproclaim
 
In contrast, modules with non-simple tops do have proper layer-stable
degenerations in general (see \cite{\preview}).

From Corollary 4.3, it is intuitively clear that absence of proper
layer-stable degenerations of the modules represented by $\grassSS$ does
not suffice to guarantee existence of a moduli space for the
representations with radical layering $\SS$; a confirmation will follow
in 4.B.  In light of this fact, the next theorem distinguishes $\toptd$
and its closed subvarieties among the quasi-affine subvarieties of
$\modlad$.  (Recall that $\R$ denotes the representation map from
$\toptd$ to the set of isomorphism classes of $d$-dimensional modules
with top $T$.)

\proclaim{Theorem 4.4}  Suppose that $T$ is squarefree and $W$ a closed
$\GL_d$-stable subvariety of the classical variety $\toptd$ of
$d$-dimensional modules with top $T$.  Let $V$ be the corresponding
closed
$\autlap$-stable subvariety of $\grasstd$ (under the bijection of
Proposition {\rm 2.5}  --  in particular, $W = \toptd$ implies $V =
\grasstd$).  Then the following statements are equivalent:
\roster
\item"(1)"  There exists a fine moduli space for the  isomorphism
classes of modules represented by $W$, i\.e., for the isomorphism
classes $R(x)$, $x \in W$.

\item"(2)" There exists a coarse moduli space for the  isomorphism
classes of modules represented by $W$. 

\item"(3)"  $V$ \underbar{is} the fine moduli space for the modules
represented by $W$. 

\item"(4)"  $W$ has a geometric quotient modulo $\GL_d$.

\item"(5)"  All $R(x)$ with $x\in W$ are devoid of top-stable
degenerations. 

\item"(6)"  All $C \in V$ are fully invariant in $P$.
\endroster 
 In case the modules represented by $W$
possess a moduli space, they all split into local direct summands.
Moreover, restriction to $V$ of the family $(\Gamma,\gamma)$ over
$\grasstd$ constructed in Lemma {\rm 3.3} is the universal
family in this situation.   

In particular:  If the the $d$-dimensional
modules with top $T$ have a moduli space classifying them up to
isomorphism, then the universal family of such modules is parametrized by
$\grasstd$ and restricts to the trivial bundle on each of the
affine patches $\grassS$.   
\endproclaim

\demo{Proof} That (3)$\implies$(1)$\implies$(2) is trivial, as is
(4)$\implies$(5).  Moreover, (6)$\implies$(4) is due to the fact that
(6) makes $V$ its own geometric quotient modulo $\autlap$; by
Proposition 2.5, (4) follows. The equivalence of (5) and (6) can be
obtained from Theorem 4.2.  That (2) implies (5) can be deduced from
Criterion 2.1 (just restrict the locally universal bundle over $\toptd$
to $W$).

(6)$\implies$(3).  Adopt (6), which means that all
$\autlap$-orbits of $V$ are reduced to points.  Then
$V$ is its own geometric quotient modulo $\autlap$ and,  using Criterion
2.6, we infer that $V$ is a coarse moduli space for the families of
modules isomorphic to $P/C$ with $C \in V$.  To see that $V$ is even a
fine moduli space, we restrict the family
$(\Gamma,\gamma)$ parametrized by $\grasstd$, constructed in Lemma 3.3,
to the subvariety $V$ and show that this latter family is universal.  So
let $(\Delta, \delta)$ be any family of $d$-dimensional modules
recruited from the $P/C$ with $C \in V$, where
$\Delta$ is a vector bundle over a variety X say.  Define a map
$\tau: X \rightarrow V \subseteq \grasstd$ by sending any point $x$ in
$X$ to the unique point $C \in \grasstd$ with the property that the
fibre of
$\Delta$ above $x$ is $\la$-isomorphic to $P/C$.   In view of the fact 
that $V$ is already  known to be a coarse moduli space, meaning that the
bijection $\alpha =
\phi|_V^{-1}$ from the set of isomorphism classes of modules $P/C$ with
$C \in V$ to the variety $V$ satisfies the conditions of 1.6$'$ in
\cite{\New, p.24}, it is now routine to check that
$\tau$ is a morphism and
$(\Delta, \delta)$ is equivalent to $\tau^*((\Gamma, \gamma)|_V)$.
Obviously,
$\tau$ is unique with this property.  

The final statements follow from Theorem 4.2 and the previous
paragraph, respectively.       
\qed \enddemo

Squarefreeness  of $T$ cannot be discarded
from the hypotheses of Theorem 4.4. To back this up, we again refer to
Example 3.7.  Namely, let
$W \subseteq \toptd$ be the
$\GL_d$-orbit representing the module $M = \la e_1 \oplus S_1$.  Then
$W$ is closed in $\toptd$, but, while $W$ has a trivial geometric
quotient modulo $\GL_d$, condition (6) of the theorem is violated.  

Theorem 4.4 shows, in particular, that the isomorphism classes of
$d$-dimensional modules with top $T$ have a fine moduli space provided
that the simple summands of $T$ do not recur as composition factors of
$JP/(\soc JP)$.  As was pointed out to the author by Crawley-Boevey, the
following special case can also be derived from
of King's work 
\cite{\Ki}:  Namely, if $T = \la e/ Je$ is simple and
$eJe = 0$, then, for any $d \in \NN$, the modules with top $T$ and fixed
class in $\text{K}_0(\la)$  --  that is, with fixed dimension vector 
--  possess a fine moduli space.  (Indeed, if the total dimension is $d$
and $\Theta:
\text{K}_0(\la) \rightarrow \ZZ$ is the $\ZZ$-linear map which sends the
class of $T$ to $-(d -1)$ and all other simples to $1$, then the modules
specified are $\Theta$-stable.)   
 
The situation where, for all $d$, the $d$-dimensional modules with fixed
simple top $T$ have a fine moduli space, can be
characterized in terms of quiver and relations of
$\la$, a justification being immediate from Theorem 4.4.    

\proclaim{Corollary 4.5}  Given a simple module $T = \la e/Je$, the
following statements are equivalent:
\roster
\item"(a)" For all $d$, the $d$-dimensional modules with top $T$ have a
coarse moduli space.

\item"(b)" For all $d$, the variety $\grasstd$ is a fine moduli space
for the $d$-dimensional modules with top $T$. 

\item"(c)" Every left ideal $C\subseteq \la e$ is of the form $C=
C\la \cap \la e$. 
  
\item"(d)" For every element $\lambda \in Je$ and every oriented cycle
$\omega \in e J e$, the product $\lambda \omega$ belongs to $\la
\lambda$. 
\endroster

\noindent Clearly, these conditions are satisfied whenever $(Je)^2 =
0$, a fortiori when $eJe = 0$. \qed
 \endproclaim

To further illustrate the richness of the representation theory in
situations where the modules with fixed top are classifiable, we mention
that arbitrary projective varieties occur (up to
isomorphism) as unions of irreducible components of moduli spaces
$\grasstd$ for families of modules with fixed dimension $d$ and fixed
top $T$:  The method of \cite{\GeomI, Theorem G} can be adapted to show
that each projective variety $V$ can be realized as a fine moduli space
of uniserial modules with a fixed sequence $\SS = \bigl(S(1), \dots,
S(d)\bigr)$ of consecutive composition factors.  This was done by 
L\. Hille, to whom the author had communicated her construction method
prior to publication (see \cite{\Hil, Example}).  In Hille's
example, the quiver of the underlying algebra $\la$ is of a form
ensuring that $\grassSS$ is a union of irreducible components of
the encompassing variety $\grasstd$ (in general, this not the case; see
\cite{\irred}).

\subhead 4.B. Remarks on quotients of $\grassSS$ 
\endsubhead

As one might suspect, there is a plethora of cases where the
$d$-dimensional modules with top
$T$ do not possess a moduli space, whereas, for each $d$-dimensional
semisimple sequence with top $T$, the modules with radical layering
$\SS$ do.  The simplest example illustrating this fact is Example 2.4: 
If
$\SS = (S_1, S_1, S_2)$ and $\SS' = (S_1, S_1 \oplus S_2, 0)$, each of 
the two subvarieties $\grassSS$ and $\grass{\SS'}$ trivially has a
moduli space (a singleton in each case), whereas
$\operatorname{\frak{Grass}}^{S_1}_3 = \grassSS \cup \grass{\SS'}$ does
not (by Theorem 4.4).  

From Criterion 2.6, we know that the representations with radical
layering
$\SS$ have a coarse moduli space precisely when the $\autlap$-space
$\grassSS$ has a categorical quotient which is an orbit map.  Obviously,
closedness of the $\autlap$-orbits in
$\grassSS$ (equivalently, absence of proper layer-stable degenerations
among the modules with radical layering $\SS$) is a necessary condition
for this event.  In view of Corollary 4.3, it is automatically satisfied
when $T$ is simple (but not so for more general
$T$; see \cite{\preview}).  Yet, even for a sequence $\SS$ with simple
top $T$,  the modules with radical layering $\SS$ frequently fail to
possess a coarse moduli space.  Indeed, another obstacle is as follows: 
Assume, for the moment that
$\grassSS$ is irreducible and $\pi: \grassSS \rightarrow \grassSS/
\autlap$ a categorical quotient which is an orbit map.  Then $\pi$ is
dominant, and hence, by 
\cite{\Bor, AG\.10\.1}, all $\autlap$-orbits of $\pi$ have the same
dimension.  For an example where this condition is violated, let $\l 
\ge 2$ and consider the algebra 
$$\la = KQ / \langle \text{all paths of length\ } l+1 \rangle,$$  where
$Q$ is the quiver

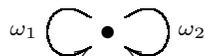
\ignore
$$\xymatrixrowsep{2pc}\xymatrixcolsep{4pc}
\xymatrix{
\bullet \ar@(ul,dl)_{\omega_1} \ar@(ur,dr)^{\omega_2} }$$
\endignore
\medskip

\noindent If $S$ is the unique simple left $\la$-module and $\SS$ the
$(l+1)$-dimensional semisimple sequence $(S, S, \dots, S)$, then
$\grassSS = (\PP^1)^l \times \AA^{l (l-1) / 2}$ is irreducible, but all
numbers between $0$ and $l-1$ arise as dimensions of
$\autlap$-orbits of $\grassSS$; for details, see
\cite{\GeomII, Example}.    

On the other hand, we believe such skips in the orbit dimension to be
the only impediment (beyond nonclosedness of the orbits for nonsimple
$T$) in the way of a moduli space.   In fact, a variety of examples led
us to the following sharper conjecture.       

\definition{Conjecture 4.6}  The following conditions are equivalent for
a semisimple sequence $\SS$ with squarefree top $T$:

(1) All $\autlap$-orbits of $\grassSS$ are closed and the orbit
dimension is constant on each irreducible component of $\grassSS$.

(2)  There exists a coarse moduli space for the representations with
radical layering $\SS$.

(3)  There exists a fine moduli space for the representations with
radical layering $\SS$.

(4)  $\grassSS$ has a geometric quotient modulo $\autlap$. 

\smallskip

If $T$ is simple, the first requirement under (1) is automatic, and the
second is equivalent to constant dimension of the endomorphism rings
$\End_\la(P/C)$, where $C$ traces an irreducible component of
$\grassSS$.
\enddefinition

\head 5. Algebras of finite local representation type   \endhead

Our main objective in this section is to characterize the algebras
which, up to isomorphism, permit only finitely many modules with a fixed
simple top $T$.  We start by considering an invariant that measures the
size of the category $\Add^T_\la$ consisting of the direct sums of left
$\la$-modules with top $T$.  This invariant fails to be left-right
symmetric in general.

\definition{Definition 5.1}  Given a simple left $\la$-module $T$ with
projective cover $P$, the {\it local dimension of\/}
$\lamod$ {\it at \/} $T$, denoted
$$\locdim^T (\lamod),$$ is the maximum of the following differences
$$\dim \C \ - \text{\ generic fibre dimension of\ } \phi \text{\ on\ }
\C,$$ where $\C$ runs through the irreducible components of the
varieties $\grasstd$ for $1 \le d \le \dim_K P$.  (Recall that
$\phi$ stands for the pertinent
 representation map.)
\enddefinition

The local dimension of $\la$ at $T$ can be obtained from the
$\autlap$-stable affine covers of the varieties $\grasstd$ for $d
\ge 1$.  Indeed, $\locdim^T(\lamod)$ equals the maximum of the
differences `$\dim \D  -  \text{generic\ } \autlap \text{-orbit
dimension on \ }\D$', where $\D$ traces the irreducible components of
the $\grassS$ and $\sigma$ ranges through all skeletons with
top $T$.  For a computation of the generic fibre dimensions on the
irreducible components
$\D$ via a (fairly small) system of linear equations, we refer to
\cite{\preview}.  Yet, in checking whether $\locdim^T(\lamod) = 0$ for
an algebra $\la$ presented in terms of quiver and relations, it is
usually more convenient to use one of the last two of the equivalent
conditions presented in Theorem 5.2 below and combine the procedures
of Section 3.B with Gr\"obner methods.  In case
$\lamod$ has vanishing local dimension at $T$, all local modules with 
top $T$ can be computed from quiver and relations of $\la$.  Due to the
theorem, they are classified by their radical layerings.  Our
concluding example demonstrates that, by contrast, dimension vectors do
not separate non-isomorphic candidates in general.       

\proclaim{Theorem 5.2} For any simple left $\la$-module
$T$ with projective cover $P$, the following statements are equivalent:
\smallskip

{\rm (1)}  There are only finitely many left
$\la$-modules with top $T$, up to isomorphism.
\smallskip

{\rm (2)} $\locdim^T (\lamod) = 0$.
\smallskip

{\rm (3)} If $M$ and $N$ are any $\la$-modules with top $T$, then
$\SS(M)=
\SS(N)$ implies $M\cong N$.
\smallskip

{\rm (4)} For every semisimple sequence $\SS$ with top $T$, the variety
$\grassSS$ is either empty or  irreducible of dimension
$$\mu_T(M) - \dim_K \End_\la (M)$$ for some (equivalently, for all) $M$
with radical layering $\SS$.
\smallskip

{\rm (5)}  For every skeleton $\S$ with top $T$, the variety $\grassS$
is either empty or consists of a  single $\autlap$-orbit.
\smallskip

If these equivalent conditions are satisfied, the varieties
$\grassSS$, where
$\SS$ ranges through the semisimple sequences with top $T$, are either
empty or isomorphic to full affine spaces. If, in addition, $Q$ does not
contain oriented cycles, then each $\grassSS$ is either empty or a point.
\endproclaim

\demo{Proof}  (1)$\implies$(2).  Adopt (1) and assume $\SS$ to be a
semisimple sequence with top $T$ such that $\grassSS \ne \varnothing$. 
If $\C$ is an irreducible component of $\grassSS$, then $\C$ consists of
only finitely many fibres of $\phi$ by hypothesis.  All of them are
closed in $\grassSS$ by Corollary 4.3, and a fortiori closed in $\C$. 
Therefore $\C$ consists of a single fibre.  In particular, $\dim \C$
equals the generic fibre dimension of $\phi$ on $\C$.  Since, for each
$d$, the irreducible components of $\grasstd$ are among the closures, in
$\grasstd$, of the irreducible components of the $\grassSS$'s covering
$\grasstd$, we conclude that $\locdim^T (\lamod) = 0$.

(2)$\implies$(1).  Given (2), it suffices to prove that, for any
semisimple sequence $\SS$ with top $T$, there are only finitely many
modules
$M$ with $\SS(M) = \SS$ up to isomorphism; indeed, there are just
finitely many such sequences $\SS$, since their total dimensions are all
bounded by $\dim_K P$.  Assume
$\grassSS
\ne
\varnothing$ and let $\C$ be an irreducible component of $\grassSS$.  If
$C \in \C$ has an $\autlap$-orbit of maximal dimension, then
$\dim \autlap.C = \dim \C$ by hypothesis, and since $\autlap.C$ is
irreducible and closed in $\grassSS$ by Proposition 4.1, we infer that
$\C = \autlap.C$, i.e.,
$\phi(\C)$ is a singleton.

(1)$\implies$(3).  We assume (1), let $\SS = (\SS_0, \dots, \SS_L)$
be any semisimple sequence with top $T$ such that $\grassSS$ is
nonempty, and $C \in \grassSS$.   We show by induction on the dimension
$d$ of $\SS$ that
$M = P/C$ is the only local module with semisimple sequence $\SS$, up to
isomorphism.  Let $k$ be the largest integer with $\SS_k \ne 0$. Without
loss of generality, $d \ge 2$, which means that $k \ge 2$, because
$\dim_K \SS_0 = \dim_K T = 1$.  Then $M' = M/J^kM$ is local with
semisimple sequence $\SS(M') = (\SS_0, \dots, \SS_{k-1}, 0,
\dots, 0)$, and by the induction hypothesis, $M'$ is unique with this
property.  This focuses our attention on the following subset of
$\modlad$, namely on
$$\Ext = \{x \in \modlad \mid  \exists \ \text{an exact sequence}\ 0
\rightarrow \SS_k \rightarrow R(x) \rightarrow M'
\rightarrow 0\}.$$ By  \cite{\BongAdv, 6.3}, $\Ext$ is an irreducible
subset of
$\modlad$.  Therefore $\Ext \cap \toptd$ is irreducible as well; indeed,
due to the upper semicontinuity of the maps $\dim_K \hom_\la( -, S_r)$,
this intersection is open in $\Ext$.  Moreover, our construction
guarantees that 
$$\operatorname{\bold{Mod}}(\SS) = \Ext \cap
\operatorname{\bold{Mod}}(\SS) =
\bigl(\Ext \cap \toptd
\bigr) \setminus
\bigl(\bigcup_{\underline{\dim}\ \SS' > \underline{\dim}\ \SS}
\operatorname{\bold{Mod}}(\SS') \bigr).$$  Since the union
$\bigcup_{\underline{\dim}\ \SS' > \underline{\dim}\ \SS}
\operatorname{\bold{Mod}}(\SS')$ is closed in $\toptd$ by Lemma 2.8(i)
and Proposition 2.5, we infer that $\operatorname{\bold{Mod}}(\SS)$ is
again irreducible.  Hence so is
$\grassSS$ by Proposition 2.5.  In light of our hypothesis, $\grassSS$
consists of finitely many orbits, one of which is $\autlap.C$, and the
fact that all of these orbits are closed in $\grassSS$ by part (3) of
Proposition 4.1 thus yields
$\grassSS = \autlap.C$.  This finishes the induction and establishes (3).

(3)$\iff$(4) is immediate from Proposition 4.1(3), as is the fact that
(3) implies the supplementary statement.   The implications
(3)$\implies$(5)$\implies$(2) are obvious. 
\qed
\enddemo

As a special case we recover the following result of Bongartz in
\cite{\Bong}: 

\proclaim{Corollary 5.3}  Let 
$\SS = (\SS_0, \dots, \SS_L)$ be a semisimple sequence with the property
that each nonzero $\SS_i$ is simple.  If $\la$ has only finitely many
uniserial modules, up to isomorphism, there is at most one uniserial left
$\la$-module with sequence $\SS$ of consecutive composition factors.\qed
\endproclaim

The following extension of Example 2.4 shows that, even in case
$\locdim^T (\lamod) = 0$, the varieties $\grasstd$ may have arbitrarily
large dimension.

\definition{Example 5.4}  For any $m \ge 1$, we present  a finite
dimensional algebra $\la$, together with a simple left
$\la$-module
$T$, such that $\locdim^T(\lamod) = 0$ while
$\operatorname{\frak{Grass}}^T_{m+2}$ contains a copy of $\PP^m$. Let
$\la = Q/I$, where $Q$ is the quiver

\ignore
$$\xymatrixrowsep{1pc}\xymatrixcolsep{3pc}
\xymatrix{ 1 \ar@(ul,dl)_{\omega} \ar[r]^{\alpha} &2 }$$
\endignore

\noindent and $I$ the ideal generated by $\omega^{m+1}$.  Moreover,
let $T =
S_1$.  Using Theorem 5.2, one readily verifies that
$\locdim^T(\lamod) = 0$.  Via the methods of \cite{\preview}, one
moreover checks that $\operatorname{\frak{Grass}}^T_{m+2}$ contains
precisely $m+1$ distinct orbits, all representing modules with dimension
vector
$(m+1,1)$; these are isomorphic to $\AA^0, \AA^1, \dots,
\AA^m$ and correspond to the submodules $C_j = \sum_{0 \le i
\le m,\ i \ne j}
\la \alpha \omega^i$ of $\la e_1$, respectively.  From Lemma 2.8 and
Theorem 4.2 one thus derives by induction that the orbit closures are
stacked into one another as follows:
$$\overline{\autlap.C_{j}} = \bigcup_{i\le j} \autlap.C_i.$$  Here is a
sketch:  Given that $\autlap.C_1 \cong \AA^1$, this orbit fails to be
closed, and hence its closure is
$\autlap.C_1 \cup \autlap.C_0$.  Similarly, the closure of $\autlap.C_2$
is seen to contain $\autlap.C_1$, and so on.  In particular, we conclude
that the closure of
$\autlap.C_m$ is isomorphic to $\PP^m$.\qed
\enddefinition

\head 6. Proof of Theorem 3.14 \endhead

We start by showing that $\chi_\S$ is well-defined and surjective and
then check that $\psi_\S$ is a bijective set map.  That $\psi_\S$ is
actually an isomorphism of varieties can be verified by an argument
modeled on that of Bongartz and the author in 
\cite{\GeomII, Theorem A}; we leave it to the reader to make the obvious
modifications.  That the triangle of maps commutes, is obvious.

It will be convenient to write $\S^{(r)}$ for the set of paths in $\S$
which start in the vertex $e_r$ and to let 
$\S_l$, resp\. $\S^{(r)}_l$, be the subset of $\S$, resp\. of $\S^{(r)}$,
consisting of the paths of length
$l$; accordingly, $\S_l = \bigcup_{1 \le r \le t} \S^{(r)}_l$ and  $\S =
\bigcup_{0 \le l \le L} \S_l$.
   
For well-definedness of $\chi_\S$, we need to show that, for any point
$c \in \AffS$, the factor module
$M(c) = P / U(c)$ has skeleton $\S$.  Set $x_r = e_r + U(c) \in M(c)$ 
for $1 \le r \le t$.  Clearly,
$M(c)$ is generated by the elements $p x_r$ for $r \le t$ and 
$p \in \S^{(r)}$.  Indeed, with the argument backing Proposition 3.12, one
checks that, given any path
$q$ and $r\le t$, the element $q x_r$ of $M(c)$
is a
$K$-linear combination of terms $p x_s$ with $p \in \S^{(s)}$, $s
\le t$.  So it suffices to verify that, for each $l \in
\{0, \dots, L\}$, the subset
$\bigcup_{1 \le r \le t} \{p x_r \mid p \in \S^{(r)}_l\}$ of
$M(c)$ is linearly independent modulo
$J^{l+1} M(c)$.  On the assumption that this fails, let $l_0$ be minimal
with respect to failure for some factor algebra of $\la$ and some point
$c$ in a suitable $\AffS$.  It is harmless to assume that
$l_0 = L$.  For otherwise we can enlarge the ideal $I
\subseteq KQ$ so as to include all paths of length $l_0 +1$,  replace
$\S$ by 
$\bigl(\bigcup_{l \le l_0} \S_l \bigr)$, and let $d$ be the cardinality
of this union; in this modified setup
$l_0$ is still minimal with respect to failure of our independence
condition.  So we only need to refute the assumption that the set
$\S^{(1)}_L x_1 \cup
\cdots \cup
\S^{(t)}_L x_t$ is linearly dependent in $M(c)$. 

To obtain a convenient framework for a comparison of coefficients,
consider the projective left ideal
$\Phat := KQ e_1 \oplus \cdots \oplus KQ e_t$ of $KQ$ and note
that $\Phat/I \Phat \cong P$.  Further, we let
$V(c)$ be the submodule of $\Phat$ generated by the differences  $\alpha
p  \ - \ \sum \Sb q\in
\S(\alpha,p)
\endSb c_{\alpha p, q} q$, where $(\alpha, p)$ ranges through the
$\S$-critical pairs.  Viewed as a
$KQ$-module, $M(c)$ is then isomorphic to the quotient $\Phat/
\bigl( V(c)+ I\Phat \bigr)$.  

Next, we note that, for $1 \le r \le t$, the set of paths in
$KQ e_r$ is the disjoint union of the following two sets:  The set
$\S^{(r)}$, and that of all paths of the form $u \alpha p$,
where $(\alpha, p)$ is a $\S$-critical pair with $p \in \S^{(r)}$ and
$u$ a path of length
$\ge 0$.  This permits us to define a $K$-linear transformation 
$$F_c : \Phat \longrightarrow \Phat$$  
as follows:  Suppose that $p$ is
any path starting in one of the vertices $e_r$ for $r \le t$.  We  set
$F_c(p) = p$ if $p \in
\S^{(r)}$; moreover, we define 
$F_c(u \alpha p)$ = 
$\sum \Sb q \in \S(\alpha,p) \endSb c_{\alpha p, q} u q$  in the second
case.  Our choice of $c$ in $\AffS$ guarantees that $F_c^L (I\Phat) =
0$.  Indeed, if $v$ is any path, of length $l$ say, Proposition 3.12 yields a congruence
$$v \seq \sum_{q \in \S}
\tau^v_{q}(X) q$$  in the polynomial ring $\A$ over $KQ$, where the
$\tau^v_{q}(X)$ are suitable polynomials in
$K[X_{\alpha p ,q}]$.  Along the induction backing Proposition 3.12, we then
find
$F_c^l(v) = \sum_{q \in \S}
\tau^v_{q}(c) q$. Hence, if $\rho \in I\Phat = Ie_1 + \cdots + Ie_t$,
then $F_c^L(\rho)= \sum_{q \in \S} \tau^\rho_{q}(c) q =0$ because the
$\tau^\rho_{q}(X)$ annihilate all points in
$\AffS$ by construction. It follows that $F_c^L (I\Phat) = 0$.  Next we
observe that, for each element
$z \in \Phat$, the difference $z - F_c(z)$ belongs to
$V(c)$, which makes $V(c)$ invariant under $F_c$.  

By assumption, there is a finite list $p_1,
\dots, p_s$ of distinct paths in $\S$, say $p_i \in \S^{(r_i)}$, 
together with nonzero scalars
$k_1, \dots, k_s$, such that the linear combination
$\sum_{1 \le i \le s} k_i p_i x_{r_i}$ vanishes in $M(c)$.  
 On moving to the above presentation of
$M(c)$ as a left
$KQ$-module, we infer that the element 
$\sum_{1 \le i \le s} k_i p_i$ of $\Phat$ is a sum
$z + a$ with $z \in V(c)$ and $a \in I \Phat$.  Now $F_c^L$ annihilates
$a$ and maps $z$ back to $V(c)$, while leaving the sum $\sum_{1 \le i
\le s} k_i p_i$ fixed, thus yielding $\sum_{1 \le i \le s} k_i p_i  \in
V(c)$.  Since any path including a non-route on
$\S$ as a right subpath is again a non-route,  this yields an equality
$$\sum_{1 \le i \le s} k_i p_i \ \ = \ \  \sum_{ ( \alpha, p
)\ \S\text{-critical}}\  
\sum_{k \ge 1}
\ b_{\alpha p, k} u_{\alpha p, k}
\biggl(\alpha p \ -
\sum \Sb q \in \S(\alpha,p) \endSb c_{\alpha p, q} q 
\biggr), \tag\dagger$$   for suitable scalars $b_{\alpha p,
k}$ and paths $u_{\alpha p, k}$ of lengths
$\ge 0$, starting in the endpoint of $\alpha$, respectively; clearly, we
may assume that $u_{\alpha p, j} \ne u_{\alpha p, k}$ for $j \ne k$. 
This being an equality  in $\Phat$, we are now in a position to compare
coefficients. 

In doing this, the concept of a route on $\S$, introduced in Remark
3.15, will come in handy.   Since none of the paths 
$u_{\alpha p, k} \alpha p$ on the right-hand side of $(\dagger)$ belongs
to $\S$,  each
$p_i$ on the left-hand side must equal one of of the paths
$u_{\alpha p, k} q$.  Moreover, one    
 observes that, whenever we have an equality $p_i = u_{\alpha p, k} q$
for some $q \in
\S(\alpha,p)$, the path $u_{\alpha p, k} \alpha p$ is a route on $\S$.  
We can  therefore find a
$\S$-critical pair
$(\alpha_0, p_0)$ such that $p_0 \in \S$ has  minimal length with
respect to the following property:  There exists an index, say $k = 1$,
with  
\roster
\item"$\bullet$" $b_{\alpha_0 p_0, 1} \ne 0$ and such that

\item"$\bullet$" $u_{\alpha_0 p_0, 1} \alpha_0 p_0$ is a route on
$\S$.
\endroster

Set $w = u_{\alpha_0 p_0, 1} \alpha_0 p_0$ and note that $w \notin \S$. 
The left-hand side of 
$(\dagger)$ being a $K$-linear combination of paths in $\S$,  
 the path $w$ must cancel out of the right-hand side.   Observe that
$w$ does not equal any path of the form  $u_{\alpha p, k}
\alpha p$ with $(\alpha,p) \ne (\alpha_0,p_0)$, for $p_0$ is the longest
right subpath of $w$ which belongs to $\S$; nor does $w$ coincide with
one of the paths 
$u_{\alpha_0 p_0, k} \alpha_0 p_0$ for $k \ne 1$.  Consequently, $w$
must be one of the $u_{\alpha p, k} q$ for some
$\S$-critical pair $(\alpha,p)$, some choice of $q \in
\S(\alpha,p)$, and some $k$ with $b_{\alpha p, k} \ne 0$.  In
particular, this makes 
$u_{\alpha p, k} q$ a route on $\S$, which entails that $u_{\alpha p, k}
\alpha p$ is also a route on
$\S$.  We infer $\len(p) \ge
\len(p_0)$ by the minimality of $\len(p_0)$, and further deduce
$\len(q) > \len(p_0)$ because $q$ belongs to
$\S(\alpha,p)$.   On the other hand, the equality 
$$w = u_{\alpha_0 p_0, 1} \alpha_0 p_0 = u_{\alpha p, k} q$$  implies
$\len(p_0) \ge \len(q)$,  once again due to the maximal length of $p_0$
as a  right subpath of $w$ that belongs to
$\S$.  This contradiction shows our original assumption $\sum_{i
\le s} k_i p_i x_{r_i}\allowmathbreak =0$ to be absurd.  Hence $\S$ is a
skeleton of
$M(c)$ as claimed, which proves well-definedness of $\chi_\S$.

To see that $\chi_\S$ is a surjection, let $M$ be any $\la$-module with
skeleton $\S$; in particular, this ensures that $M$ is
$d$-dimensional with top $T$.  Choose a sequence  $y_1, \dots, y_t$ of
top elements of $M$ with $y_r = e_r y_r$.  For each
$l$, the products $p y_r$ for $p \in \S^{(r)}_l$ and $1 \le r \le t$
then form a basis for $J^lM / J^{l+1}M$.  Moreover, for any
$\S$-critical pair $(\alpha,p)$, with $p \in \S^{(r)}$ say, let
$c_{\alpha p,q}$ be the unique scalars such that
$\alpha p\, y_r = \sum_{s \le t} \sum_{q \in \S(\alpha,p)\cap \S^{(s)}}
c_{\alpha p,q} q \,y_s$.  It clearly suffices to show that the
corresponding point $c = (c_{\alpha p, q})
\in \AA^N$ belongs to the variety $\AffS$, for it is then clear that
$M \cong P / U(c)$; just send $y_r$ to $e_r + U(c)$.  So let us check $c
\in \AffS$.  For that purpose,  suppose $v \in I e_r$ for some
$r$ between $1$ and $t$.  Viewing $v$ as an element of $\A$, consider the
expansion
$v  \seq \sum_{q \in \S} \tau^{v}_{q}(X)\, q$ with
$\tau^{v}_{q}(X) \in K[X]$  supplied by Proposition 3.12.  In tandem with the
inductive procedure justifying the lemma, one obtains the following
equalities in $M$: 
$$v\, y_r = \sum_{s \le t}
\sum_{q \in \S^{(s)}} \tau^{v}_{q}(c)\, q\, y_s;$$ here we refer to the
induced $KQ$-structure of $M$.  Consequently, the linear independence of
the elements $q\, y_s$, for
$q \in \S^{(s)}$ and $s \le t$, yields vanishing of all the values
$\tau^{v}_{q}(c)$. This means that $c \in \AffS$ as asserted and shows
that $\chi_{\sigma}$ is indeed a surjection.

That the map $\psi_{\sigma}$ is injective, follows from the first part of
the proof;  indeed, it suffices to observe that, for any $c \in \AffS$,
the paths in
$\sigma$ are linearly independent modulo $U(c)$.  To verify
surjectivity, let $C \in \grassS$.  Then $P/C$ has skeleton $\S$ and top
elements $e_r + C$.  Therefore, as explained in the surjectivity
argument for
$\psi_{\S}$, there exists a point $c \in \AffS$ with the property that,
for any $\S$-critical pair $(\alpha,p)$ with $p \in \S^{(r)}$, 
$$\alpha p (e_r + C) = \sum_{s \le t} \sum_{q \in \S(\alpha,p) \cap
\S^{(s)}} c_{\alpha p,q} q (e_s + C).$$ In other words, $U(c) \subseteq
C$.  But we already know that
$U(c)$ has codimension $d$ in $P$ as does $C$, and hence $\psi_\S(c)=
U(c)= C$. This completes the proof of Theorem 3.14. \qed   

\head 7.  Moduli spaces for representations with arbitrary top $T$ -- a
preview
\endhead

Suppose that $T = \bigoplus_{1 \le i \le n} S_i^{t_i}$ is an arbitrary 
semisimple module, with projective cover $P$.  Recall from
Observation 2.3 that 
$$\autlap \cong \U \rtimes \aut_\la(T) \cong 
\U \rtimes \prod_{1 \le i \le t} \GL_{t_i}(K),$$ 
where  $\U$ denotes
the unipotent radical of $\autlap$.  Moreover, let $\T$ be the
direct product of the tori of diagonal matrices in the $\GL_{t_i}(K)$.
To realize $\T$ as a subgroup of $\autlap$, write $P = \bigoplus_{1 \le i \le n} \bigoplus_{1 \le j \le
t_i} \la x_{ij}$ for suitable elements $x_{ij} = e_i x_{ij}$.  Then $\T$
can be identified with the automorphisms of $P$ sending $x_{ij}$ to
$a_{ij} x_{ij}$, where
$\bigl(a_{ij}\bigr)_{i \le n,\,j \le t_i} \in (K^*)^{\dim T}$.  We
observe that the equality 
$\autlap = \U \rtimes \T$ holds precisely when $T$ is squarefree; it is
this fact which singles out the squarefree case in the context of the
moduli problem.  We only give a rough picture of the general case,
addressing solvability of the moduli problem for all $d$, to illustrate
the strong pressure exerted by this demand on those indecomposable
projective $\la$-modules whose radical factors occur with multiplicity
$>1$ in $T$. 

\proclaim{7.1} The following conditions are equivalent for the
semisimple module $T$ containing the $S_i$ with multiplicity $t_i$:

\roster
\item For each $d$, all $\autlap$-orbits of $\grasstd$ are closed. 

\item For each $d$, all $\U \rtimes \T$-orbits of $\grasstd$ are
singletons, and every indecomposable projective module $\la e_i$ with
$t_i > 1$ is uniserial.

\item For each $d$ and each point $C \in \grasstd$, the module $P/C$ is
a direct sum of local modules, $\dim \hom_\la (P, JP/C) = \dim \hom_\la
(P/C, JP/C)$, and all $\la e_i$ with $t_i > 1$ are uniserial.
\endroster
\endproclaim

This restricts `global' classifiability of the representations with
top $T$ to a rather narrow selection of algebras $\la$.  On the other
hand, if one splits up the class of top-$T$ modules in terms of
an invariant somewhat  reminiscent of the `words' underlying
classification of the representations of biserial algebras, one
obtains a fine moduli space for each of the resulting subclasses under
much more general circumstances.

By a slight abuse of language, we use the term {\it distinguished
basis\/} of a module
$P/C$ for any family  $\sigma = \bigl( \sigma_{ij} \bigr)_{i
\le n,\, j\le t_i}$ of sets of paths in $KQ$ with the property that
each
$\sigma_{ij}$ consists of paths starting in the vertex
$e_i$ and is closed under right subpaths (in particular, this entails
$e_i \in \sigma_{ij}$) and such that, moreover, the residue classes
$px_{ij} + C$ for $p \in \sigma_{ij}$ (with $i,j$ tracing all eligible
choices) form a
$K$-basis of
$P/C$.  Clearly, each module $P/C$ with $C \subseteq JP$ has at least
one such distinguished basis, and there are only finitely many 
families $\sigma$ that qualify as distinguished bases. 
Finally, when $\sum_{i,j} |\sigma_{ij}| = d$, let $\C(\sigma)$ be the
set of all points $C$ in $\grasstd$ such that $P/C$ has distinguished
basis
$\sigma$ (we refer to $\sigma$ as $d$-dimensional in that case).  Note
that, for each choice of $\sigma$, the set $\C(\sigma)$ is an open
subvariety of $\grasstd$, being the intersection of a
standard open affine subset of the classical Grassmannian with
$\grasstd$. (However, the $\C(\sigma)$ fail to be stable under the $\U
\rtimes \T$-action and, a fortiori, under the $\autlap$-action in
general.)

\proclaim{7.2} The following statements are equivalent for given $T$ and
$d$:

\roster
\item For each $d$-dimensional distinguished basis $\sigma$, the
family of top-$T$ modules $P/C$ with $C \in \C(\sigma)$ has a fine
moduli space. 

\item All $\U \rtimes \T$-orbits of $\grasstd$ are singletons.

\item For each $C \in \grasstd$, the quotient $P/C$ is a direct
sum of local modules and $\dim \hom_\la (P, JP/C) = \dim \hom_\la
(P/C, JP/C)$. 
\endroster
\endproclaim

If the equivalent conditions of 7.2 are satisfied,  
each of the varieties $\C(\sigma)$ is the fine moduli space
for the corresponding class $\bigl( P/C \bigr)_{C \in \C(\sigma)}$.  
Moreover,
$\C(\sigma)$ is a direct product of subvarieties of
Grassmannians of modules with simple tops,
which brings the problem back, full circle, to the simplest squarefree
case.

\enddocument